\documentclass[submission]{FPSAC2022}

\newtheorem{prop}{Proposition}
\newtheorem{thm}{Theorem}

\newtheorem{remark}{Remark}

\usepackage{lipsum}

\usepackage{tikz}
\usetikzlibrary{shapes,shadows,arrows,positioning,graphs}
\usetikzlibrary{calc,decorations.pathmorphing,intersections}
\usetikzlibrary{arrows.meta}

\title[PC Prographs and Triangulations of the Sphere]{Product-Coproduct Prographs and Triangulations of the Sphere}

\author[Nicolas Borie \and Justine Falque]{Nicolas Borie \and Justine Falque}

\address{LIGM, Univ Gustave Eiffel, CNRS, ESIEE Paris, Marne-la-Vallee, France}

\received{\today}

\abstract{In this paper, we explain how the classical Catalan families of objects  
involving paths, tableaux, triangulations, parentheses configurations and more
generalize canonically to a three-dimensional version. In 
particular, we present product-coproduct prographs as central objects 
explaining the combinatorics of the triangulations of the sphere. Then we expose
a natural way to extend the Tamari lattice to the product-coproduct prographs.}

\resume{Dans ce papier, nous expliquons comment les familles d'objets 
combinatoires comptées par les nombres de Catalan, comme les chemins, tableaux, 
triangulations, parenthésages, ou autre, se généralisent 
dans le monde des objets comptés par les nombres de Catalan tridimensionnels. En
particulier, nous présentons les prographes produit-coproduit comme centraux pour
expliquer la combinatoire des triangulations de la sph\`ere. Nous exposons alors une 
manière naturelle d'étendre le treillis de Tamari aux prographes produit-coproduit.} 


\usepackage[backend=bibtex]{biblatex}
\begin{document}

\maketitle

\section{Introduction}

This paper is about a generalization of the
famous combinatorial problem defining Catalan numbers and involving 
parentheses configurations, binary trees and triangulations of the $n$-gon. 
This paper is a sequel of~\cite{borie_fpsac_slc_2017} where 
product-coproduct prographs (PC prographs) were introduced
to study the combinatorics of the three-dimensional Catalan world.
The work presented in this paper follows the results obtained in 
FPSAC 2017~\cite{borie_fpsac_slc_2017}, whose statements can be 
summarized in Table~\ref{tab_first_prop}.

In this paper, we complete these first statements to the dual side of planar objects 
(triangulations of the sphere) and the combinatorics of the planar objects. 

Mireille Bousquet-M\'elou did expose in~\cite{MBM_planar_maps} that the three-dimensional
Catalan numbers count the number of triangulations of the bipolar sphere. By embedding 
PC prographs on the sphere and linking the only global output to the global input, we
present PC prographs as Voronoï diagrams of the triangulations of the bipolar sphere. 
This implies that the original bijection from rectangular Young Tableaux to PC prographs
can be extended to triangulations.

\begin{center}
  \begin{footnotesize}
    \begin{tabular}{|cc|} \hline

    \multicolumn{2}{|c|}{\textbf{Integer sequences and formal power series}} \\ 
    Catalan numbers & Three-dimensional Catalan numbers \\ 
    $\displaystyle\frac{(2n)!}{(n+1)! \cdot n!}$ & 
    $\displaystyle\frac{2 \cdot (3n)!}{n! \cdot (n+1)! \cdot (n+2)!}$ \\ 

    \multicolumn{2}{|c|}{\textbf{First combinatorial classes of objects}} \\ 
    \qquad \qquad Two-rows standard Young tableaux \qquad \qquad & Three-rows standard Young tableaux \\ 
    $$
    \begin{tabular}{|c|c|c|c|} \hline
      & $\cdots$ & $\cdots$ & $2n$ \\ \hline
      $1$ & $\cdots$ & $\cdots$ & \\ \hline
    \end{tabular}
    $$
    &
    $$ 
    \begin{tabular}{|c|c|c|c|} \hline
      $\quad$ & $\cdots$ & $\cdots$ & $3n$ \\ \hline
      & $\cdots$ & $\cdots$ &  \\ \hline
      $1$ & $\cdots$ & $\cdots$ & \\ \hline
      \multicolumn{4}{c}{$\quad$} \\
    \end{tabular}
    $$
    \\ 

    \multicolumn{2}{|c|}{\textbf{Realization as Operators rules in Algebras}} \\ 
    Parentheses configurations of an $n$-product & Ways to assemble-disassemble $n$-times each \\
    in an associative algebra & in an associative and coassociative bialgebra \\ 
    & $\times \cdot \Delta \cdot \times \cdot \Delta$ \\ 
    $( \bullet \times \bullet ) \times \bullet$ & $\times \cdot (\times \otimes Id) \cdot (Id \otimes \Delta) \cdot \Delta$ \\ 
     & $\times \cdot (Id \otimes \times) \cdot (Id \otimes \Delta) \cdot \Delta$ \\ 
    $\bullet \times ( \bullet \times \bullet )$ & $\times \cdot (\times \otimes Id) \cdot (\Delta \otimes Id) \cdot \Delta$ \\
    & $\times \cdot (Id \otimes \times) \cdot (\Delta \otimes Id) \cdot \Delta$ \\ 

    \multicolumn{2}{|c|}{\textbf{Planar combinatorial objects}} \\
    Binary trees & PC prographs\\ 
    \begin{tikzpicture}[scale=0.55]
      \tikzstyle{prod}=[fill,draw,rectangle,minimum size=5pt,inner sep=1pt]
      \tikzstyle{cop}=[fill,draw,circle,minimum size=6pt,inner sep=1pt]
    
      \draw (16,0.5) -- (16,0);
      \draw (16,0) -- (14,-2);
      \draw (16,0) -- (18,-2);
      \draw (17,-1) -- (16,-2);
      \draw[white] (16,-2) -- (16,-3);
    
      \draw (16,0) node[prod] (p2) {$~$};
      \draw (17,-1) node[prod] (p1) {$~$};
    \end{tikzpicture}
    & 
    \begin{tikzpicture}[scale=0.55]
      \tikzstyle{prod}=[fill,draw,rectangle,minimum size=5pt,inner sep=1pt]
      \tikzstyle{cop}=[fill,draw,circle,minimum size=6pt,inner sep=1pt]
    
      \draw (16,0.5) -- (16,0);
      \draw (16,0) -- (15,-1) -- (17, -3);
      \draw (16,0) -- (18,-2) -- (17, -3);
      \draw (17,-1) -- (16,-2);
      \draw (17,-3) -- (17,-3.5);
    
      \draw (16,0) node[prod] (p2) {$~$};
      \draw (17,-1) node[prod] (p1) {$~$};
      \draw (16,-2) node[cop] (c2) {$~$};
      \draw (17,-3) node[cop] (c1) {$~$};
    
    \end{tikzpicture}
    \\ \hline

    \end{tabular}
  \end{footnotesize}
  \vspace{0.2cm}
  \textbf{Table 1:} Classical Catalan world and its three-dimensional counterpart.~\label{tab_first_prop}
\end{center}

By duality, since rotations in prographs correspond to triangle flipping,
we define four rotations rules on PC-prographs. We explain that these rotations 
correspond exactly to all admissible triangulation flips on the triangulations side.
Among these four rotations rules, two are extentions of the classical 
rotations balancing binary search trees and the other two are suggested by 
the geometry of the triangulations of the sphere.
Finally, we present some structural properties of the set of PC prographs of size $n$
extended with the four rotations viewed as oriented rewriting rules. We argue that, 
in some cases, PC prographs can be viewed as a gluing of two binary trees; endowing
this subfamily with two of the four rotations makes it a lattice that is a product of
the Tamari lattice by itself.
However, although adding the other two rotations does generate all PC-prographs
with faithful action of the rotations, this does not provide a lattice.






This paper is organised as follows. In the next Section, we first recall what 
a PC prograph is and describe the
duality between PC prographs and rooted triangulations of the sphere. We 
exploit the duality to introduce and orient flips and rotations. In 
Section~\ref{structures}, we investigate structures that can be defined over the 
set of PC prographs.

\section{Duality and triangulations}

PC prographs are connected oriented planar graphs without cycles whose vertices 
are products or coproducts. Coproducts have one input and two
outputs and products have two inputs and one output. Just one
coproduct (resp.\ product) has its input (resp.\ output) unconnected, which is 
the input (resp.\ output) of the prograph (see Figure~\ref{Fig1} for examples).
Since each vertex has degree three, the dual graph consists of triangles.

Now, embed a prograph on the sphere and connect its output to
its input going through the exterior non visible side of the sphere. The
prograph becomes a Voronoï diagram of some triangulation of the sphere.
Since~\cite{MBM_planar_maps, bipolar} prove that bipolar
oriented triangulations of the sphere are counted by the three-dimensional
Catalan numbers, this procedure provides an alternative and explicit bijection.

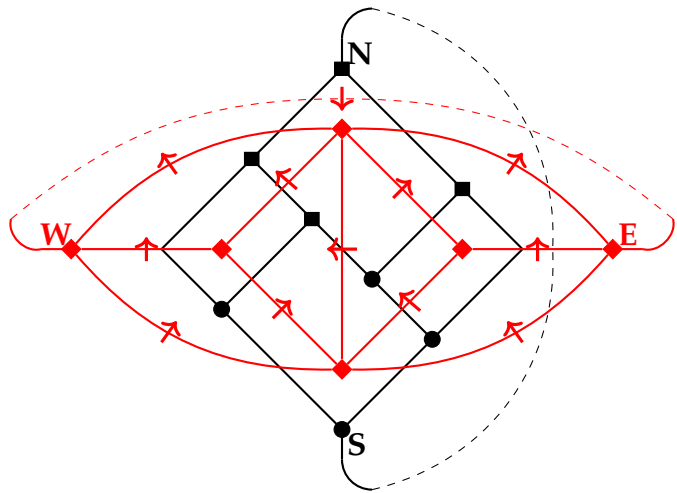
\begin{figure}
  \centering
  \begin{tikzpicture}[scale=0.4]
    \tikzstyle{prod}=[fill,draw,rectangle,minimum size=5pt,inner sep=1pt]
    \tikzstyle{cop}=[fill,draw,circle,minimum size=6pt,inner sep=1pt]
    \tikzstyle{vor}=[red,fill=red,draw,diamond,minimum size=6pt,inner sep=1pt]

    \draw (0,0) node[cop] (c1) {$~$};
    \draw (-4,4) node[cop] (c2) {$~$};
    \draw (3,3) node[cop] (c3) {$~$};
    \draw (1,5) node[cop] (c4) {$~$};

    \draw (-1,7) node[prod] (p1) {$~$};
    \draw (-3,9) node[prod] (p2) {$~$};
    \draw (4,8) node[prod] (p3) {$~$};
    \draw (0,12) node[prod] (p4) {$~$};
    
    \draw[thick] (0, -1) -- (0, 0) -- (6, 6) -- (0, 12) -- (0, 13);
    \draw[thick] (0, 0) -- (-6, 6) -- (0, 12);
    \draw[thick] (3, 3) -- (-3, 9);
    \draw[thick] (-4, 4) -- (-1, 7);
    \draw[thick] (1, 5) -- (4, 8);

    \draw (0, 13) edge[thick, bend left=45] (1, 14);
    \draw (1, 14) edge[dashed, bend left=40] (7, 6);
    \draw (0, -1) edge[thick, bend right=45] (1, -2);
    \draw (1, -2) edge[dashed, bend right=40] (7, 6);

    \draw (0,2) node[vor] (v1) {$~$};    
    \draw (0,10) node[vor] (v2) {$~$};
    \draw (-4,6) node[vor] (v3) {$~$};
    \draw (4,6) node[vor] (v4) {$~$};
    \draw (-9,6) node[vor] (v5) {$~$};
    \draw (9,6) node[vor] (v6) {$~$};

    \draw (v1) edge[thick, red] (v3);
    \draw (v1) edge[thick, red] (v4);
    \draw (v2) edge[thick, red] (v3);
    \draw (v2) edge[thick, red] (v4);
    \draw (v1) edge[thick, red] (v2);

    \draw (v3) edge[thick, red] (v5);
    \draw (v1) edge[thick, red, bend left=25] (v5);
    \draw (v2) edge[thick, red, bend right=25] (v5);
    \draw (v4) edge[thick, red] (v6);
    \draw (v2) edge[thick, red, bend left=25] (v6);
    \draw (v1) edge[thick, red, bend right=25] (v6);

    \draw (9, 6) edge[thick, red] (10, 6);
    \draw (10, 6) edge[thick, red, bend right=60] (11, 7);
    \draw (11, 7) edge[dashed, red, bend right=20] (0, 11);
    \draw (-9, 6) edge[thick, red] (-10, 6);
    \draw (-10, 6) edge[thick, red, bend left=60] (-11, 7);
    \draw (-11, 7) edge[dashed, red, bend left=20] (0, 11);

    \draw (0.5, 6) edge[->, very thick, red] (-0.5, 6);
    \draw (2.6, 3.95) edge[->, very thick, red] (1.9, 4.55);
    \draw (-1.5, 8.0) edge[->, very thick, red] (-2.2, 8.6);
    \draw (-2.3, 3.7) edge[->, very thick, red] (-1.7, 4.3);
    \draw (1.7, 7.7) edge[->, very thick, red] (2.3, 8.3);
    \draw (6.5, 5.6) edge[->, very thick, red] (6.5, 6.4);
    \draw (-6.5, 5.6) edge[->, very thick, red] (-6.5, 6.4);

    \draw (5.5, 8.5) edge[->, very thick, red] (6, 9.2);
    \draw (-5.5, 8.5) edge[->, very thick, red] (-6, 9.2);    
    \draw (6, 2.9) edge[->, very thick, red] (5.5, 3.6);
    \draw (-6, 2.9) edge[->, very thick, red] (-5.5, 3.6);
    \draw (0, 11.4) edge[->, very thick, red] (0, 10.6);

    \draw (-9.5, 6.5) node[red] (north) {$\mathbf{W}$};
    \draw (9.5, 6.5) node[red] (south) {$\mathbf{E}$};
    \draw (0.5, -0.5) node (east) {$\mathbf{S}$};
    \draw (0.6, 12.5) node (west) {$\mathbf{N}$};
    
  \end{tikzpicture}
  \caption{A self-dual prograph (black) and its bipolar-oriented triangulation (red).\label{duality_exemple}}
\end{figure}

\begin{figure}[h]
  \centering
\begin{tikzpicture}[scale=0.55]
  \tikzstyle{prod}=[fill,draw,rectangle,minimum size=4pt,inner sep=1pt]
  \tikzstyle{cop}=[fill,draw,circle,minimum size=2pt,inner sep=1pt]
  \tikzstyle{vor}=[fill=red,red,draw,diamond,minimum size=2pt,inner sep=1pt]

  \draw (0, 0.5) -- (0, 0);
  \draw (0,0) -- (-0.5,-0.5) -- (0, -1);
  \draw (0,0) -- (0.5,-0.5) -- (0, -1);
  \draw (0, -1) -- (0, -2);
  \draw (0,-2) -- (-0.5,-2.5) -- (0, -3);
  \draw (0,-2) -- (0.5,-2.5) -- (0, -3);
  \draw (0, -3) -- (0, -3.5);

  \draw (0,0) node[prod] (p2) {$~$};
  \draw (0,-1) node[cop] (c2) {$~$};
  \draw (0,-2) node[prod] (p1) {$~$};
  \draw (0,-3) node[cop] (c1) {$~$};

  \draw (0, -0.5) node[vor] (v1) {$~$};
  \draw (0, -2.5) node[vor] (v2) {$~$};
  \draw (-1, -1.5) node[vor] (v3) {$~$};
  \draw (1, -1.5) node[vor] (v4) {$~$};
  \draw[red] (-1.5, -1.5) -- (1.5, -1.5);
  \draw[red] (-1, -1.5) -- (0, -0.5) -- (1, -1.5);
  \draw[red] (-1, -1.5) -- (0, -2.5) -- (1, -1.5);

  \draw (6,0.5) -- (6,0);
  \draw (5,-1) -- (6,0) -- (7, -1);
  \draw (4,-2) -- (5,-1) -- (6, -2);
  \draw (7,-1) -- (6,-2);
  \draw (4,-2) -- (5,-3) -- (6,-2);
  \draw (5,-3) -- (5,-3.5);

  \draw (6,0) node[prod] (p2) {$~$};
  \draw (5,-1) node[prod] (p1) {$~$};
  \draw (6,-2) node[cop] (c2) {$~$};
  \draw (5,-3) node[cop] (c1) {$~$};

  \draw (5, -2) node[vor] (v1) {$~$};
  \draw (6, -1) node[vor] (v2) {$~$};
  \draw (3.5, -1.5) node[vor] (v3) {$~$};
  \draw (7.5, -1.5) node[vor] (v4) {$~$};
  \draw[red] (3, -1.5) -- (3.5, -1.5);
  \draw[red] (7.5, -1.5) -- (8, -1.5);
  \draw (v3) edge[red] (v1);
  \draw (v1) edge[red] (v2);
  \draw (v2) edge[red] (v4);
  \draw (v3) edge[red, bend left=40] (v2);
  \draw (v1) edge[red, bend right=40] (v4);

  \draw (12,0.5) -- (12,0);
  \draw (12,0) -- (13.5,-1.5) -- (12, -3);
  \draw (12,0) -- (10.5,-1.5) -- (12, -3);
  \draw (11,-1) -- (11.5,-1.5) -- (11, -2);
  \draw (12,-3) -- (12,-3.5);

  \draw (12,0) node[prod] (p2) {$~$};
  \draw (11,-1) node[prod] (p1) {$~$};
  \draw (11,-2) node[cop] (c2) {$~$};
  \draw (12,-3) node[cop] (c1) {$~$};

  \draw (11, -1.5) node[vor] (v1) {$~$};
  \draw (10, -1.5) node[vor] (v2) {$~$};
  \draw (12, -1.5) node[vor] (v3) {$~$};
  \draw (14, -1.5) node[vor] (v4) {$~$};
  \draw[red] (9.5, -1.5) -- (14.5, -1.5);
  \draw (10, -1.5) edge[red, bend left=40] (11, -0.5);
  \draw (11, -0.5) edge[red, bend left=40] (12, -1.5);
  \draw (10, -1.5) edge[red, bend right=40] (11, -2.5);
  \draw (11, -2.5) edge[red, bend right=40] (12, -1.5);

  \draw (3,-5.5) -- (3,-6);
  \draw (3,-6) -- (4.5,-7.5) -- (3, -9);
  \draw (3,-6) -- (1.5,-7.5) -- (3, -9);
  \draw (4,-7) -- (3.5,-7.5) -- (4, -8);
  \draw (3,-9) -- (3,-9.5);

  \draw (3,-6) node[prod] (p2) {$~$};
  \draw (4,-7) node[prod] (p1) {$~$};
  \draw (4,-8) node[cop] (c2) {$~$};
  \draw (3,-9) node[cop] (c1) {$~$};

  \draw (1, -7.5) node[vor] (v1) {$~$};
  \draw (3, -7.5) node[vor] (v2) {$~$};
  \draw (4, -7.5) node[vor] (v3) {$~$};
  \draw (5, -7.5) node[vor] (v4) {$~$};
  \draw[red] (0.5, -7.5) -- (5.5, -7.5);
  \draw (3, -7.5) edge[red, bend left=40] (4, -6.5);
  \draw (4, -6.5) edge[red, bend left=40] (5, -7.5);
  \draw (3, -7.5) edge[red, bend right=40] (4, -8.5);
  \draw (4, -8.5) edge[red, bend right=40] (5, -7.5);

  \draw (9,-5.5) -- (9,-6);
  \draw (9,-6) -- (8,-7) -- (10, -9);
  \draw (9,-6) -- (11,-8) -- (10, -9);
  \draw (10,-7) -- (9,-8);
  \draw (10,-9) -- (10,-9.5);

  \draw (9,-6) node[prod] (p2) {$~$};
  \draw (10,-7) node[prod] (p1) {$~$};
  \draw (9,-8) node[cop] (c2) {$~$};
  \draw (10,-9) node[cop] (c1) {$~$};

  \draw (10, -8) node[vor] (v1) {$~$};
  \draw (9, -7) node[vor] (v2) {$~$};
  \draw (7.5, -7.5) node[vor] (v3) {$~$};
  \draw (11.5, -7.5) node[vor] (v4) {$~$};
  \draw[red] (7, -7.5) -- (7.5, -7.5);
  \draw[red] (11.5, -7.5) -- (12, -7.5);
  \draw (v3) edge[red] (v2);
  \draw (v1) edge[red] (v2);
  \draw (v1) edge[red] (v4);
  \draw (v3) edge[red, bend right=40] (v1);
  \draw (v2) edge[red, bend left=40] (v4);

\end{tikzpicture}
\caption{The five prographs with two coproducts, two products (in black),
  and their associated triangulations (in red).}~\label{Fig1}
\end{figure}
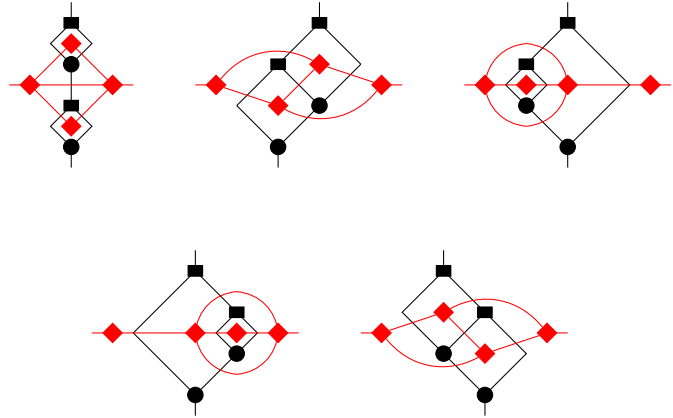

For example, consider the following three-row and rectangular standard Young tableau:
\begin{displaymath}
  \begin{array}{|c|c|c|c|} \hline
    7 & 8 & 11 & 12 \\ \hline
    3 & 4 & 9 & 10 \\ \hline
    1 & 2 & 5 & 6 \\ \hline
  \end{array}.
\end{displaymath}
This tableau is stable by the Schützenberger
involution~(see~\cite{borie_fpsac_slc_2017, falque_hard_bij}) 
(rotate $180^{\circ}$ the tableau then complement its values 
sending $i$ to $3n+1 - i$). Using the central bijection of~\cite{borie_fpsac_slc_2017} 
from tableaux to prographs, we obtain the black prograph of
Figure~\ref{duality_exemple}, with coproducts represented as circles and
products as squares. Considering this prograph as a Voronoï diagram, we
build the red oriented triangulation of the sphere with a single edge
passing through its dark side. This red dashed edge (its midpoint could
have been sent to infinity) splits the dark side of the sphere into two
triangles: one containing the North pole of the sphere and the other one
containing its South pole.

\begin{remark}
  Switching the context from prographs to bipolar oriented
  triangulations of the sphere, the Schützenberger involution
  consists in rotating the sphere by 180 degrees, swapping the poles,
  and flipping the orientation of the edges of the triangulation. It is also
  equivalent to the antipodal map of the sphere with renamed
  poles since they have been swapped.
\end{remark}

The reader can check that the triangulation obtained in
Figure~\ref{duality_exemple} is stable by the Schützenberger
involution.

\begin{thm}
  With PC prographs as an intermediate object, graph duality provides an 
  explicit and constructive bijection from 3-row rectangular standard Young
  tableaux onto bipolar-oriented triangulations of the
  sphere. This bijection commutes with the Schützenberger involution.
\end{thm}

Since there exists numerous algorithms for exhaustive enumeration
and (possibly random) generation of Young tableaux, this bijection
using PC prographs enables a lot of calculations
and computations on the triangulations side. For example, we believe that
the easiest way to check if two oriented triangulations of the bipolar
sphere are homotopic is to build their associated prographs with graph
duality and check whether they are associated with the same Young tableau.

\begin{remark}
The world of prographs meets here the world of oriented planar maps. 
Indeed, prographs are planar assemblies of
operators with identified inputs and outputs; if one puts reasonable operators on
the vertices of a planar oriented map, with respect to the number of ongoing
and outgoing edges for each vertex, one builds a prograph. Since, for PC
prographs, a single entry means the vertex is a coproduct and two entries mean
it is a product, there is no choice to be made and everyone works on the same objects.
\end{remark}

\section{Flips and rotations}

The Tamari lattice has numerous different realizations\,: using parentheses 
configurations, binary trees, ordered forests, triangulations
and more. Our first attempts to uncover a structure on PC prographs were by
focusing on the prographs and Young tableaux side. Unfortunately, there were
too many choices when trying to establish a set of transforms and
give them orientations, thereby exhibiting a partial order (or more).

On the other hand, the geometric approach revealed the needed rigidity. 
When one considers triangulations, the reasonable flip consists in focusing
on an edge, check that
it is the border of two triangles; if the union of these two triangles
forms a quadrilateral and the edge is a diagonal of it, then one flips this
edge by removing it and drawing the other diagonal instead. Using this
rule alone and
by exhaustion of cases, we now show that there is only one reasonable
way to consider flips and rotations.

By graph duality, each triangle contains a single operator, a product or
a coproduct. Pairs of neighboring triangles (sharing at least an edge) can be 
enumerated by the entries inside the associated Young tableau
since, if triangles are neighbors, there exists an edge between both 
involed operators. We forget $1$ which
virtually connects the northernmost product with
the southernmost coproduct \emph{via} the dark side of the sphere. This
pointed edge in PC prographs correspond to a pointed edge in triangulations
that we will not flip (it plays the same role as the first side of the
$n$-gon in the classical Catalan world and Tamari lattice.)

In PC prographs, there are three types of outputs\,: outputs of products,
left outputs of coproducts and right outputs of coproducts. Similarly, there
are three types of inputs\,: inputs of coproducts, left inputs of products
and right inputs of products. Since edges in PC prographs connect
an output of operator to an input, we should consider nine cases.

\begin{figure}
\begin{tabular}{|c|c|c|c|} \hline
  type of        & input                  & left input             & right input           \\ 
  edges          & of coproduct           & of product             & of product            \\ \hline
  
  left output    & Type I                 & Type II                & Type III              \\ 
  of coproduct   & \textbf{reduced}       & no possible flip       & reducible to type VII \\ \hline
  
  right output   & Type IV                & Type V                 & Type VI               \\ 
  of coproduct   & reducible to type I    & \textbf{reduced}       & no possible flip      \\ \hline
  
  output         & Type VII               & Type VIII              & Type IX               \\ 
  of product     & reducible to type V    & reducible to type IX   & \textbf{reduced}      \\ \hline
\end{tabular}
\caption{An edge is flippable depending on the type of its entry and its outgoing.}~\label{tabflip}
\end{figure}

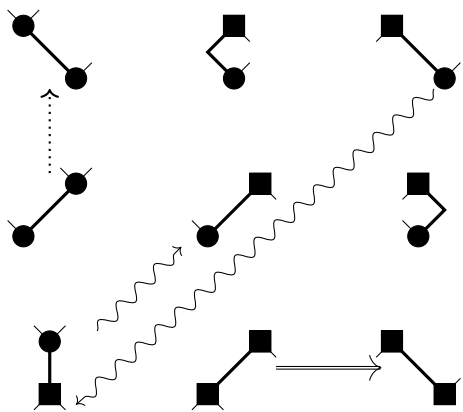
\begin{figure}
    \centering
    \begin{tikzpicture}[scale=0.7]
      \tikzstyle{prod}=[fill,draw,rectangle,minimum size=8pt,inner sep=1pt]
      \tikzstyle{cop}=[fill,draw,circle,minimum size=8pt,inner sep=1pt]

      \draw (0,-1) node[cop] {};
      \draw[very thick] (0, -1) -- (-1, 0);
      \draw (-1,0) node[cop] {};
      \draw (0, -1) -- (0.3, -0.7);
      \draw (-1, 0) -- (-1.3, 0.3);
      \draw (-1, 0) -- (-0.7, 0.3);

      \draw (3,-1) node[cop] {};
      \draw[very thick] (3, -1) -- (2.5, -0.5) -- (3, 0);
      \draw (3, 0) node[prod] {};
      \draw (3, -1) -- (3.3, -0.7);
      \draw (3, 0) -- (3.3, -0.3);

      \draw (7,-1) node[cop] {};
      \draw[very thick] (7, -1) -- (6, 0);
      \draw (6,0) node[prod] {};
      \draw (7, -1) -- (7.3, -0.7);
      \draw (6, 0) -- (5.7, -0.3);

      \draw (-1,-4) node[cop] {};
      \draw[very thick] (-1, -4) -- (0, -3);
      \draw (0,-3) node[cop] {};
      \draw (-1, -4) -- (-1.3, -3.7);
      \draw (0, -3) -- (-0.3, -2.7);
      \draw (0, -3) -- (0.3, -2.7);

      \draw (2.5,-4) node[cop] {};
      \draw[very thick] (2.5, -4) -- (3.5, -3);
      \draw (3.5, -3) node[prod] {};
      \draw (2.5, -4) -- (2.2, -3.7);
      \draw (3.5, -3) -- (3.8, -3.3);

      \draw (6.5, -4) node[cop] {};
      \draw[very thick] (6.5, -4) -- (7, -3.5) -- (6.5, -3);
      \draw (6.5, -3) node[prod] {};
      \draw (6.5, -4) -- (6.2, -3.7);
      \draw (6.5, -3) -- (6.2, -3.3);

      \draw (-0.5,-7) node[prod] {};
      \draw[very thick] (-0.5, -7) -- (-0.5, -6);
      \draw (-0.5,-6) node[cop] {};
      \draw (-0.5, -7) -- (-0.8, -7.3);
      \draw (-0.5, -7) -- (-0.2, -7.3);
      \draw (-0.5, -6) -- (-0.8, -5.7);
      \draw (-0.5, -6) -- (-0.2, -5.7);

      \draw (2.5,-7) node[prod] {};
      \draw[very thick] (2.5, -7) -- (3.5, -6);
      \draw (3.5, -6) node[prod] {};
      \draw (2.5, -7) -- (2.2, -7.3);
      \draw (2.5, -7) -- (2.8, -7.3);
      \draw (3.5, -6) -- (3.8, -6.3);

      \draw (7, -7) node[prod] {};
      \draw[very thick] (7, -7) -- (6, -6);
      \draw (6, -6) node[prod] {};
      \draw (7, -7) -- (7.3, -6.7);
      \draw (6, -6) -- (5.7, -6.3);

      \draw[->, thick, dotted] (-0.5, -2.8) -- (-0.5, -1.2);

      \draw[->, double] (3.8, -6.5) -- (5.8, -6.5);

      \draw[->, decorate, decoration={snake}] (0.4,-5.8) -- (2, -4.2);

      \draw[->, decorate, decoration={snake}, bend right=50] (6.8, -1.2) -- (0, -7.2);

    \end{tikzpicture}
  \caption{The 9 types of edges and their reductions via oriented rotations}~\label{edgeflip}
  \end{figure}

The only choice we made is the orientation of operations. 
Here, types I, V
and IX are reduced by choice, and, since we want to be coherent with the
Schützenberger involution, the only other possible choice is to reverse
all orientations simultaneously. Reversing all rules, types III, IV and VI would
have been the reduced edges. This corresponds to exchanging the left
and right everywhere and flipping the east and west on the sphere.

The reader can check exhaustively that a flip of types II and VI would create a
loop inside the prograph, whatever the choice of orientation.
Therefore flipping these two types of edges is technically impossible. 
It so happens that these edges correspond exactly
to pathological organizations of triangles on the sphere. For example, on
the sphere, there exist valid triangulations in which two neighboring triangles 
have two edges in common. The union of two such triangles does not form 
a quadrilateral and no flip is possible.

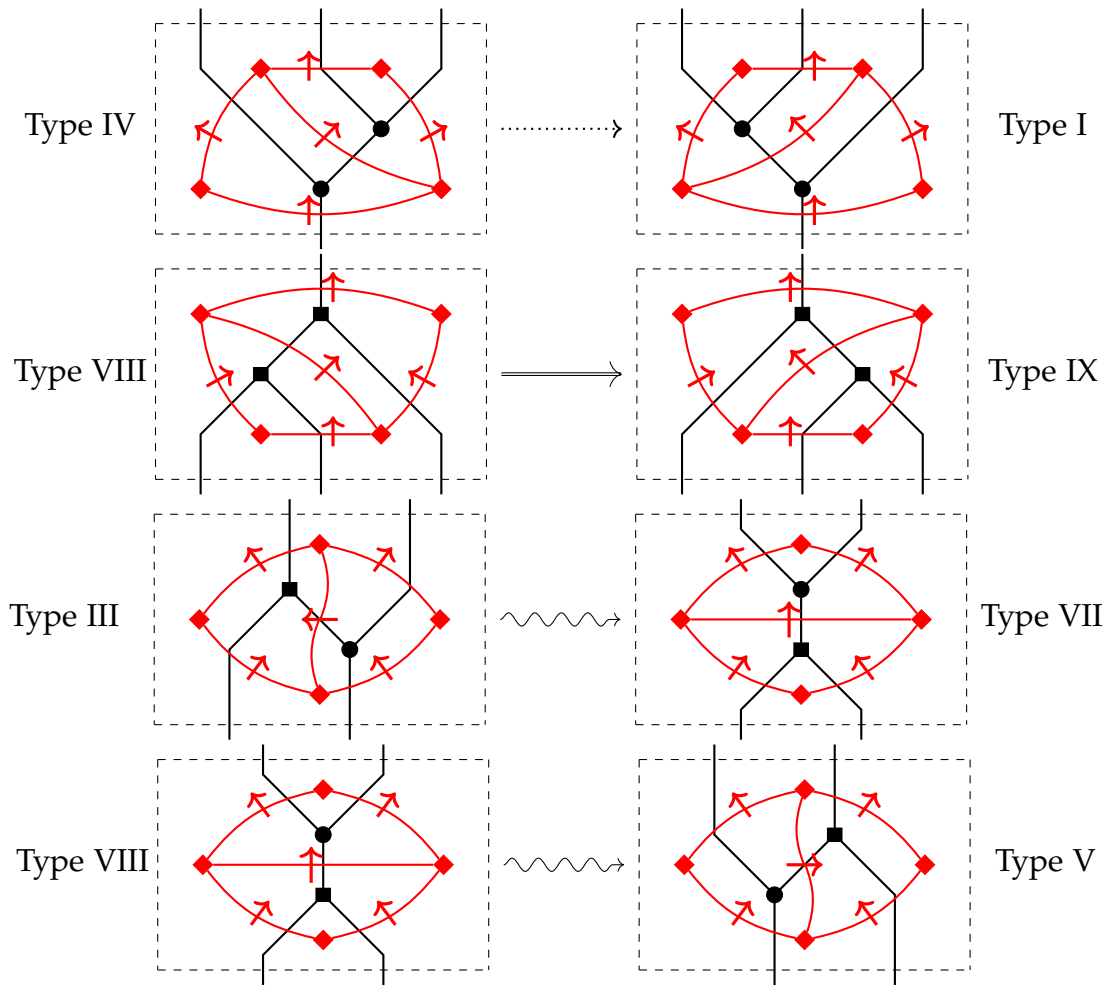
\begin{figure}[h]
  \centering
  \begin{tikzpicture}[scale=0.8]
    \tikzstyle{prod}=[fill,draw,rectangle,minimum size=5pt,inner sep=1pt]
    \tikzstyle{cop}=[fill,draw,circle,minimum size=6pt,inner sep=1pt]
    \tikzstyle{vor}=[red,fill=red,draw,diamond,minimum size=6pt,inner sep=1pt]

    \draw (-2,1) node (t4) {Type IV};
    \draw (14,1) node (t1) {Type I};
    
    \draw (3,1) node[cop] (p2) {$~$};
    \draw (2,0) node[cop] (p1) {$~$};

    \draw[thick] (2, -1) -- (2, 0) -- (4, 2) -- (4, 3);
    \draw[thick] (3, 1) -- (2, 2) -- (2, 3);
    \draw[thick] (2, 0) -- (0, 2) -- (0, 3);

    \draw (1,2) node[vor] (v1) {$~$};
    \draw (3,2) node[vor] (v2) {$~$};
    \draw (0,0) node[vor] (v3) {$~$};
    \draw (4,0) node[vor] (v4) {$~$};

    \draw (v1) edge[thick, red] (v2);
    \draw (v1) edge[thick, red, bend right=20] (v3);
    \draw (v2) edge[thick, red, bend left=20] (v4);
    \draw (v3) edge[thick, red, bend right=20] (v4);
    \draw (v1) edge[thick, red, bend right=20] (v4);

    \draw[dashed] (-0.75,2.75) -- (4.75, 2.75);
    \draw[dashed] (-0.75,-0.75) -- (4.75, -0.75);
    \draw[dashed] (-0.75,2.75) -- (-0.75, -0.75);
    \draw[dashed] (4.75,2.75) -- (4.75, -0.75);

    \draw (1.9,0.7) edge[->, very thick, red] (2.3, 1.1);
    \draw (1.8, -0.6) edge[->, very thick, red] (1.8, -0.1);
    \draw (1.8, 1.8) edge[->, very thick, red] (1.8, 2.3);
    \draw (0.35, 0.8) edge[->, very thick, red] (-0.1, 1.05);
    \draw (3.65, 0.8) edge[->, very thick, red] (4.1, 1.05);
    
    \draw[->, thick, dotted] (5,1) -- (7,1);
    
    \draw (9,1) node[cop] (p2) {$~$};
    \draw (10,0) node[cop] (p1) {$~$};

    \draw[thick] (10, -1) -- (10, 0) -- (12, 2) -- (12, 3);
    \draw[thick] (9, 1) -- (10, 2) -- (10, 3);
    \draw[thick] (10, 0) -- (8, 2) -- (8, 3);

    \draw (9,2) node[vor] (v1) {$~$};
    \draw (11,2) node[vor] (v2) {$~$};
    \draw (8,0) node[vor] (v3) {$~$};
    \draw (12,0) node[vor] (v4) {$~$};

    \draw (v1) edge[thick, red] (v2);
    \draw (v1) edge[thick, red, bend right=20] (v3);
    \draw (v2) edge[thick, red, bend left=20] (v4);
    \draw (v3) edge[thick, red, bend right=20] (v4);
    \draw (v2) edge[thick, red, bend left=20] (v3);

    \draw[dashed] (7.25,2.75) -- (12.75, 2.75);
    \draw[dashed] (7.25,-0.75) -- (12.75, -0.75);
    \draw[dashed] (7.25,2.75) -- (7.25, -0.75);
    \draw[dashed] (12.75,2.75) -- (12.75, -0.75);

    \draw (10.2,0.8) edge[->, very thick, red] (9.8, 1.2);
    \draw (10.2, -0.6) edge[->, very thick, red] (10.2, -0.1);
    \draw (10.2, 1.8) edge[->, very thick, red] (10.2, 2.3);
    \draw (8.35, 0.8) edge[->, very thick, red] (7.9, 1.05);
    \draw (11.65, 0.8) edge[->, very thick, red] (12.1, 1.05);
    
  \end{tikzpicture}

  \begin{tikzpicture}[scale=0.8]
    \tikzstyle{prod}=[fill,draw,rectangle,minimum size=5pt,inner sep=1pt]
    \tikzstyle{cop}=[fill,draw,circle,minimum size=6pt,inner sep=1pt]
    \tikzstyle{vor}=[red,fill=red,draw,diamond,minimum size=6pt,inner sep=1pt]

    \draw (-2,1) node (t8) {Type VIII};
    \draw (14,1) node (t9) {Type IX};
    
    \draw (1,1) node[prod] (p2) {$~$};
    \draw (2,2) node[prod] (p1) {$~$};

    \draw[thick] (2, 3) -- (2, 2) -- (1, 1) -- (0, 0) -- (0, -1);
    \draw[thick] (1, 1) -- (2, 0) -- (2, -1);
    \draw[thick] (2, 2) -- (4, 0) -- (4, -1);

    \draw (1,0) node[vor] (v1) {$~$};
    \draw (3,0) node[vor] (v2) {$~$};
    \draw (0,2) node[vor] (v3) {$~$};
    \draw (4,2) node[vor] (v4) {$~$};

    \draw (v1) edge[thick, red] (v2);
    \draw (v1) edge[thick, red, bend left=20] (v3);
    \draw (v2) edge[thick, red, bend right=20] (v4);
    \draw (v3) edge[thick, red, bend left=20] (v4);
    \draw (v2) edge[thick, red, bend right=20] (v3);

    \draw[dashed] (-0.75,2.75) -- (4.75, 2.75);
    \draw[dashed] (-0.75,-0.75) -- (4.75, -0.75);
    \draw[dashed] (-0.75,2.75) -- (-0.75, -0.75);
    \draw[dashed] (4.75,2.75) -- (4.75, -0.75);

    \draw (1.9, 0.9) edge[->, very thick, red] (2.3, 1.3);
    \draw (2.2, -0.2) edge[->, very thick, red] (2.2, 0.3);
    \draw (2.2, 2.2) edge[->, very thick, red] (2.2, 2.7);
    \draw (0.1, 0.8) edge[->, very thick, red] (0.55, 1.05);
    \draw (3.9, 0.8) edge[->, very thick, red] (3.45, 1.05);

    \draw[->, double] (5,1) -- (7,1);
    
    \draw (11,1) node[prod] (ap2) {$~$};
    \draw (10,2) node[prod] (ap1) {$~$};

    \draw[thick] (10, 3) -- (10, 2) -- (9, 1) -- (8, 0) -- (8, -1);
    \draw[thick] (11, 1) -- (10, 0) -- (10, -1);
    \draw[thick] (10, 2) -- (12, 0) -- (12, -1);

    \draw (9,0) node[vor] (av1) {$~$};
    \draw (11,0) node[vor] (av2) {$~$};
    \draw (8,2) node[vor] (av3) {$~$};
    \draw (12,2) node[vor] (av4) {$~$};

    \draw (av1) edge[thick, red] (av2);
    \draw (av1) edge[thick, red, bend left=20] (av3);
    \draw (av2) edge[thick, red, bend right=20] (av4);
    \draw (av3) edge[thick, red, bend left=20] (av4);
    \draw (av1) edge[thick, red, bend left=20] (av4);

    \draw[dashed] (7.25,2.75) -- (12.75, 2.75);
    \draw[dashed] (7.25,-0.75) -- (12.75, -0.75);
    \draw[dashed] (7.25,2.75) -- (7.25, -0.75);
    \draw[dashed] (12.75,2.75) -- (12.75, -0.75);

    \draw (10.2,1) edge[->, very thick, red] (9.8, 1.4);
    \draw (9.8, -0.2) edge[->, very thick, red] (9.8, 0.3);
    \draw (9.8, 2.2) edge[->, very thick, red] (9.8, 2.7);
    \draw (8.1, 0.8) edge[->, very thick, red] (8.55, 1.05);
    \draw (11.9, 0.8) edge[->, very thick, red] (11.45, 1.05);
    
  \end{tikzpicture}

  \begin{tikzpicture}[scale=0.8]
    \tikzstyle{prod}=[fill,draw,rectangle,minimum size=5pt,inner sep=1pt]
    \tikzstyle{cop}=[fill,draw,circle,minimum size=6pt,inner sep=1pt]
    \tikzstyle{vor}=[red,fill=red,draw,diamond,minimum size=6pt,inner sep=1pt]

    \draw (-2.25,1) node (t3) {Type III};
    \draw (14,1) node (t7) {Type VII};
    
    \draw (1.5,1.5) node[prod] (p2) {$~$};
    \draw (2.5,0.5) node[cop] (p1) {$~$};

    \draw[thick] (1.5, 3) -- (1.5, 1.5) -- (2.5, 0.5) -- (2.5, -1);
    \draw[thick] (1.5,1.5) -- (0.5, 0.5) -- (0.5, -1);
    \draw[thick] (2.5,0.5) -- (3.5, 1.5) -- (3.5, 3);

    \draw (0,1) node[vor] (v1) {$~$};
    \draw (4,1) node[vor] (v2) {$~$};
    \draw (2,2.25) node[vor] (v3) {$~$};
    \draw (2,-0.25) node[vor] (v4) {$~$};
    \draw (1.925,0.8) node (centeru) {};
    \draw (2.075,1.2) node (centerd) {};

    \draw (v1) edge[thick, red, bend left=20] (v3);
    \draw (v3) edge[thick, red, bend left=20] (v2);
    \draw (v2) edge[thick, red, bend left=20] (v4);
    \draw (v4) edge[thick, red, bend left=20] (v1);
    \draw (v3) edge[thick, red, bend left=20] (centeru);
    \draw (centerd) edge[thick, red, bend right=20] (v4);

    \draw[dashed] (-0.75,2.75) -- (4.75, 2.75);
    \draw[dashed] (-0.75,-0.75) -- (4.75, -0.75);
    \draw[dashed] (-0.75,2.75) -- (-0.75, -0.75);
    \draw[dashed] (4.75,2.75) -- (4.75, -0.75);

    \draw (2.3, 1) edge[->, very thick, red] (1.7, 1);
    \draw (0.8, 0) edge[->, very thick, red] (1.1, 0.4);
    \draw (3.2, 0) edge[->, very thick, red] (2.9, 0.4);
    \draw (1.1, 1.8) edge[->, very thick, red] (0.8, 2.2);
    \draw (2.9, 1.8) edge[->, very thick, red] (3.2, 2.2);
    
    \draw[->, decorate, decoration={snake}] (5,1) -- (7,1);

    \draw (10,0.5) node[prod] (p2) {$~$};
    \draw (10,1.5) node[cop] (p1) {$~$};

    \draw[thick] (9, -1) -- (9, -0.5) -- (10, 0.5) -- (11, -0.5) -- (11, -1);
    \draw[thick] (9, 3) -- (9, 2.5) -- (10,1.5) -- (11, 2.5) -- (11, 3);
    \draw[thick] (10,1.5) -- (10,0.5);
    
    \draw (8,1) node[vor] (v1) {$~$};
    \draw (12,1) node[vor] (v2) {$~$};
    \draw (10,2.25) node[vor] (v3) {$~$};
    \draw (10,-0.25) node[vor] (v4) {$~$};

    \draw (v1) edge[thick, red, bend left=20] (v3);
    \draw (v3) edge[thick, red, bend left=20] (v2);
    \draw (v2) edge[thick, red, bend left=20] (v4);
    \draw (v4) edge[thick, red, bend left=20] (v1);
    \draw (v1) edge[thick, red] (v2);
    
    \draw[dashed] (7.25,2.75) -- (12.75, 2.75);
    \draw[dashed] (7.25,-0.75) -- (12.75, -0.75);
    \draw[dashed] (7.25,2.75) -- (7.25, -0.75);
    \draw[dashed] (12.75,2.75) -- (12.75, -0.75);

    \draw (9.8, 0.7) edge[->, very thick, red] (9.8, 1.3);
    \draw (8.8, 0) edge[->, very thick, red] (9.1, 0.4);
    \draw (11.2, 0) edge[->, very thick, red] (10.9, 0.4);
    \draw (9.1, 1.8) edge[->, very thick, red] (8.8, 2.2);
    \draw (10.9, 1.8) edge[->, very thick, red] (11.2, 2.2);
    
  \end{tikzpicture}

  \begin{tikzpicture}[scale=0.8]
    \tikzstyle{prod}=[fill,draw,rectangle,minimum size=5pt,inner sep=1pt]
    \tikzstyle{cop}=[fill,draw,circle,minimum size=6pt,inner sep=1pt]
    \tikzstyle{vor}=[red,fill=red,draw,diamond,minimum size=6pt,inner sep=1pt]

    \draw (-2,1) node (t8) {Type VIII};
    \draw (14,1) node (t9) {Type V};
    
    \draw (2,0.5) node[prod] (p2) {$~$};
    \draw (2,1.5) node[cop] (p1) {$~$};

    \draw[thick] (1, -1) -- (1, -0.5) -- (2, 0.5) -- (3, -0.5) -- (3, -1);
    \draw[thick] (1, 3) -- (1, 2.5) -- (2,1.5) -- (3, 2.5) -- (3, 3);
    \draw[thick] (2,1.5) -- (2,0.5);
    
    \draw (0,1) node[vor] (v1) {$~$};
    \draw (4,1) node[vor] (v2) {$~$};
    \draw (2,2.25) node[vor] (v3) {$~$};
    \draw (2,-0.25) node[vor] (v4) {$~$};

    \draw (v1) edge[thick, red, bend left=20] (v3);
    \draw (v3) edge[thick, red, bend left=20] (v2);
    \draw (v2) edge[thick, red, bend left=20] (v4);
    \draw (v4) edge[thick, red, bend left=20] (v1);
    \draw (v1) edge[thick, red] (v2);    

    \draw[dashed] (-0.75,2.75) -- (4.75, 2.75);
    \draw[dashed] (-0.75,-0.75) -- (4.75, -0.75);
    \draw[dashed] (-0.75,2.75) -- (-0.75, -0.75);
    \draw[dashed] (4.75,2.75) -- (4.75, -0.75);

    \draw (1.8, 0.7) edge[->, very thick, red] (1.8, 1.3);
    \draw (0.8, 0) edge[->, very thick, red] (1.1, 0.4);
    \draw (3.2, 0) edge[->, very thick, red] (2.9, 0.4);
    \draw (1.1, 1.8) edge[->, very thick, red] (0.8, 2.2);
    \draw (2.9, 1.8) edge[->, very thick, red] (3.2, 2.2);
    
    \draw[->, decorate, decoration={snake}] (5,1) -- (7,1);

    \draw (10.5,1.5) node[prod] (p2) {$~$};
    \draw (9.5,0.5) node[cop] (p1) {$~$};

    \draw[thick] (10.5, 3) -- (10.5, 1.5) -- (9.5, 0.5) -- (9.5, -1);
    \draw[thick] (9.5,0.5) -- (8.5, 1.5) -- (8.5, 3);
    \draw[thick] (10.5,1.5) -- (11.5, 0.5) -- (11.5, -1);

    \draw (8,1) node[vor] (v1) {$~$};
    \draw (12,1) node[vor] (v2) {$~$};
    \draw (10,2.25) node[vor] (v3) {$~$};
    \draw (10,-0.25) node[vor] (v4) {$~$};
    \draw (10.075,0.8) node (centeru) {};
    \draw (9.925,1.2) node (centerd) {};

    \draw (v1) edge[thick, red, bend left=20] (v3);
    \draw (v3) edge[thick, red, bend left=20] (v2);
    \draw (v2) edge[thick, red, bend left=20] (v4);
    \draw (v4) edge[thick, red, bend left=20] (v1);
    \draw (v3) edge[thick, red, bend right=20] (centeru);
    \draw (centerd) edge[thick, red, bend left=20] (v4);

    \draw[dashed] (7.25,2.75) -- (12.75, 2.75);
    \draw[dashed] (7.25,-0.75) -- (12.75, -0.75);
    \draw[dashed] (7.25,2.75) -- (7.25, -0.75);
    \draw[dashed] (12.75,2.75) -- (12.75, -0.75);

    \draw (9.7, 1) edge[->, very thick, red] (10.3, 1);
    \draw (8.8, 0) edge[->, very thick, red] (9.1, 0.4);
    \draw (11.2, 0) edge[->, very thick, red] (10.9, 0.4);
    \draw (9.1, 1.8) edge[->, very thick, red] (8.8, 2.2);
    \draw (10.9, 1.8) edge[->, very thick, red] (11.2, 2.2);
    
  \end{tikzpicture}
  \caption{The four local rotations over PC prographs shown with their dual.}~\label{rota}
\end{figure}

\begin{thm}
The four rotation rules in PC prographs shown in Figure~\ref{rota} 
correspond to all valid edge flippings in rooted oriented 
triangulations of the sphere.
\end{thm}

The reader can check that the first two rotations are directly derived from 
the classical rotation rule in binary trees.

The four rotation rules, without their orientations, are completely
rigid as far as geometry is concerned. Anyway, the reader could 
argue that choices of orientation were made. Since we 
still want to be compatible with the Schützenberger involution, the choice 
of minimum and maximum must be reversed for product and coproduct.
But exactly as the classical Tamari lattice implies a choice (one
can set either the left comb tree or the right comb tree to be the maximum), 
there is an equivalent choice with PC prographs. This will result in the same 
kind of symmetries for the classical Tamari lattice and our new structure
on PC prographs.

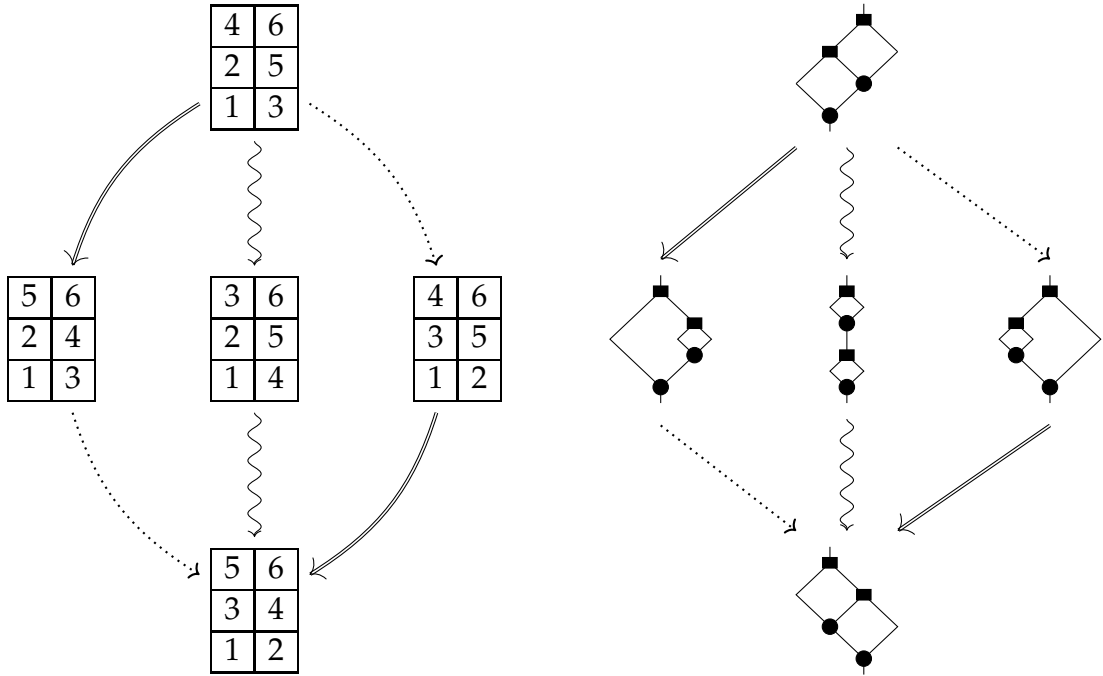
\begin{figure}[h]
  \centering
  \begin{tikzpicture}[xscale=0.9, yscale=0.85]
    \tikzstyle{prod}=[fill,draw,rectangle,minimum size=4pt,inner sep=1pt]
    \tikzstyle{cop}=[fill,draw,circle,minimum size=2pt,inner sep=1pt]

    \draw (-3, 4.25) node (up) {$ \begin{array}{|c|c|} \hline 4 & 6 \\ \hline 2 & 5 \\ \hline 1 & 3 \\ \hline  \end{array}$};

    \draw (-6, 0) node (m1) {$ \begin{array}{|c|c|} \hline 5 & 6 \\ \hline 2 & 4 \\ \hline 1 & 3 \\ \hline  \end{array}$};
    \draw (-3, 0) node (m2) {$ \begin{array}{|c|c|} \hline 3 & 6 \\ \hline 2 & 5 \\ \hline 1 & 4 \\ \hline  \end{array}$};
    \draw (0, 0) node (m3) {$ \begin{array}{|c|c|} \hline 4 & 6 \\ \hline 3 & 5 \\ \hline 1 & 2 \\ \hline  \end{array}$};

    \draw (-3, -4.25) node (dn) {$ \begin{array}{|c|c|} \hline 5 & 6 \\ \hline 3 & 4 \\ \hline 1 & 2 \\ \hline  \end{array}$};

    \draw (up) edge[->, decorate, decoration={snake}] (m2);
    \draw (up) edge[->, double, bend right=20] (m1);
    \draw (up) edge[->, thick, dotted, bend left=20] (m3);
    \draw (m2) edge[->, decorate, decoration={snake}] (dn);
    \draw (m1) edge[->, thick, dotted, bend right=20] (dn);
    \draw (m3) edge[->, double, bend left=20] (dn);

    \draw (6, 5.25) -- (6, 5);
    \draw (5.5, 4.5) -- (6, 5) -- (6.5, 4.5);
    \draw (5, 4) -- (5.5, 4.5) -- (6, 4);
    \draw (6.5, 4.5) -- (6, 4);
    \draw (5, 4) -- (5.5, 3.5) -- (6, 4);
    \draw (5.5, 3.5) -- (5.5, 3.25);
    \draw (6, 5) node[prod] (p2) {$~$};
    \draw (5.5, 4.5) node[prod] (p1) {$~$};
    \draw (6, 4) node[cop] (c2) {$~$};
    \draw (5.5, 3.5) node[cop] (c1) {$~$};
  
    \draw (5.75, 1) -- (5.75, 0.75);
    \draw (5.75,0.75) -- (5.5,0.5) -- (5.75, 0.25);
    \draw (5.75,0.75) -- (6,0.5) -- (5.75, 0.25);
    \draw (5.75, 0.25) -- (5.75, -0.25);
    \draw (5.75, -0.25) -- (5.5,-0.5) -- (5.75, -0.75);
    \draw (5.75, -0.25) -- (6,-0.5) -- (5.75, -0.75);
    \draw (5.75, -0.75) -- (5.75, -1);
    \draw (5.75, 0.75) node[prod] (p2) {$~$};
    \draw (5.75, 0.25) node[cop] (c2) {$~$};
    \draw (5.75, -0.25) node[prod] (p1) {$~$};
    \draw (5.75, -0.75) node[cop] (c1) {$~$};
    
    \draw (5.5,-3.25) -- (5.5,-3.5);
    \draw (5.5,-3.5) -- (5,-4) -- (6, -5);
    \draw (5.5,-3.5) -- (6.5,-4.5) -- (6, -5);
    \draw (6,-4) -- (5.5,-4.5);
    \draw (6,-5) -- (6,-5.25);  
    \draw (5.5,-3.5) node[prod] (p2) {$~$};
    \draw (6,-4) node[prod] (p1) {$~$};
    \draw (5.5,-4.5) node[cop] (c2) {$~$};
    \draw (6,-5) node[cop] (c1) {$~$};
  
    \draw (3, 1) -- (3, 0.75);
    \draw (3, 0.75) -- (3.75, 0) -- (3, -0.75);
    \draw (3, 0.75) -- (2.25, 0) -- (3, -0.75);
    \draw (3.5, 0.25) -- (3.25, 0) -- (3.5, -0.25);
    \draw (3, -0.75) -- (3,-1);
    \draw (3, 0.75) node[prod] (p2) {$~$};
    \draw (3.5, 0.25) node[prod] (p1) {$~$};
    \draw (3.5,-0.25) node[cop] (c2) {$~$};
    \draw (3,-0.75) node[cop] (c1) {$~$};
    
    \draw (8.75, 1) -- (8.75, 0.75);
    \draw (8.75, 0.75) -- (9.5, 0) -- (8.75, -0.75);
    \draw (8.75, 0.75) -- (8, 0) -- (8.75, -0.75);
    \draw (8.25, 0.25) -- (8.5, 0) -- (8.25, -0.25);
    \draw (8.75, -0.75) -- (8.75, -1);  
    \draw (8.75, 0.75) node[prod] (p2) {$~$};
    \draw (8.25, 0.25) node[prod] (p1) {$~$};
    \draw (8.25, -0.25) node[cop] (c2) {$~$};
    \draw (8.75, -0.75) node[cop] (c1) {$~$};

    \draw[->, decorate, decoration={snake}] (5.75, 3) -- (5.75, 1.25);
    \draw[->, double, bend right=20] (5, 3) -- (3, 1.25);
    \draw[->, thick, dotted, bend left=20] (6.5, 3) -- (8.75, 1.25);
    \draw[->, decorate, decoration={snake}] (5.75, -1.25) -- (5.75, -3);
    \draw[->, thick, dotted, bend right=20] (3, -1.35) -- (5, -3);
    \draw[->, double, bend left=20] (8.75, -1.35) -- (6.5, -3);
  \end{tikzpicture}
  \caption{The poset (lattice for this size) of PC prographs of size $2$.}~\label{assos2tab}
\end{figure}

\section{Poset structure over PC prographs}~\label{structures}

In this section, we investigate what type of structure can be placed
over the set of PC prographs.

\subsection{Gluing two trees at their canopies}~\label{glutree}

The orientation in rooted triangulations of
the sphere come directly from orientations of the classical case, the triangulations
of the rooted $n$-gon. Choosing a root side in the $n$-gon amounts to set an
orientation to a chosen side (ongoing or outgoing, the flow being
perpendicular to this side) and set the reverse orientation to all other sides.
So there are only two isomorphic choices. On the one hand, one can set the root
edge to be outgoing. The other sides are thus ongoing, hence inputs, 
and the dual of the triangulation ends up being a rooted tree of products. 
On the other hand, one can set the root edge to be 
ongoing, all other sides become outputs and the dual becomes a rooted tree 
of coproducts.

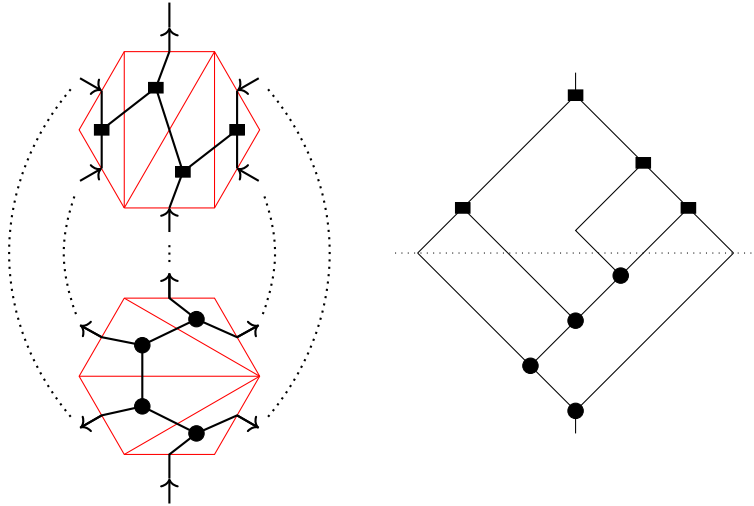
\begin{figure}
  \centering
  \begin{tikzpicture}[scale=1.2]
    \tikzstyle{prod}=[fill,draw,rectangle,minimum size=4pt,inner sep=1pt]
    \tikzstyle{cop}=[fill,draw,circle,minimum size=2pt,inner sep=1pt]

    \draw[red] (0, 1) -- (1, 1) -- (1.5, 1.866) -- (1, 2.733) -- (0, 2.733) -- (-0.5, 1.866) -- (0, 1);
    \draw[red] (0, 0) -- (1, 0) -- (1.5, -0.866) -- (1, -1.733) -- (0, -1.733) -- (-0.5, -0.866) -- (0, 0);

    \draw[red] (0, 1) -- (0, 2.733);
    \draw[red] (0, 1) -- (1, 2.733);
    \draw[red] (1, 1) -- (1, 2.733);

    \draw[red] (0, 0) -- (1.5, -0.866);
    \draw[red] (-0.5, -0.866) -- (1.5, -0.866);
    \draw[red] (0, -1.733) -- (1.5, -0.866);

    \draw (0.8, -1.5) node[cop] () {$~$};
    \draw (0.2, -1.2) node[cop] () {$~$};
    \draw (0.2, -0.522) node[cop] () {$~$};
    \draw (0.8, -0.233) node[cop] () {$~$};

    \draw[thick] (0.5, -2) -- (0.5, -1.733) -- (0.8, -1.5);
    \draw[thick] (1.49, -1.438) -- (1.25, -1.3) -- (0.8, -1.5);
    \draw[thick] (-0.49, -1.438) -- (-0.25, -1.3) -- (0.2, -1.2) -- (0.8, -1.5);
    \draw[thick] (-0.49, -0.3) -- (-0.25, -0.433) -- (0.2, -0.522) -- (0.2, -1.2);
    \draw[thick] (0.5, 0.277) -- (0.5, 0) -- (0.8, -0.233);
    \draw[thick] (1.49, -0.3) -- (1.25, -0.433) -- (0.8, -0.233) -- (0.2, -0.522);

    \draw[thick, ->] (0.5, -2.277) -- (0.5, -2);
    \draw[thick, ->] (1.25, -1.3) -- (1.49, -1.438);
    \draw[thick, ->] (-0.25, -1.3) -- (-0.49, -1.438);
    \draw[thick, ->] (-0.25, -0.433) -- (-0.49, -0.3);
    \draw[thick, ->] (0.5, 0) -- (0.5, 0.277);
    \draw[thick, ->] (1.25, -0.433) -- (1.49, -0.3);

    \draw (0.65, 1.4) node[prod] () {$~$};
    \draw (0.35, 2.333) node[prod] () {$~$};
    \draw (-0.25, 1.866) node[prod] () {$~$};
    \draw (1.25, 1.866) node[prod] () {$~$};

    \draw[thick, ->] (0.5, 0.733) -- (0.5, 1);
    \draw[thick, ->] (-0.49, 1.3) -- (-0.25, 1.438);
    \draw[thick, ->] (1.49, 1.3) -- (1.25, 1.438);
    \draw[thick, ->] (1.49, 2.438) -- (1.25, 2.3);
    \draw[thick, ->] (-0.49, 2.438) -- (-0.25, 2.3);
    \draw[thick, ->] (0.5, 2.733) -- (0.5, 3);
    \draw[thick] (0.5, 3) -- (0.5, 3.277);

    \draw[thick] (0.5, 1) -- (0.65, 1.4) -- (1.25, 1.866) -- (1.25, 1.438);
    \draw[thick] (1.25, 2.3) -- (1.25, 1.866);
    \draw[thick] (-0.25, 1.438) -- (-0.25, 2.3);
    \draw[thick] (-0.25, 1.866) -- (0.35, 2.333) -- (0.5, 2.733);
    \draw[thick] (0.65, 1.4) -- (0.35, 2.333);

    \draw (-0.49, 2.438) node (u1) {$~$};
    \draw (-0.49, 1.3) node (u2) {$~$};
    \draw (0.5, 0.733) node (u3) {$~$};
    \draw (1.49, 1.3) node (u4) {$~$};
    \draw (1.49, 2.438) node (u5) {$~$};

    \draw (-0.49, -1.438) node (d1) {$~$};
    \draw (-0.49, -0.3) node (d2) {$~$};
    \draw (0.5, 0.277) node (d3) {$~$};
    \draw (1.49, -0.3) node (d4) {$~$};
    \draw (1.49, -1.438) node (d5) {$~$};

    \draw (d1) edge[thick, dotted, bend left=40] (u1);
    \draw (d2) edge[thick, dotted, bend left=20] (u2);
    \draw (d3) edge[thick, dotted] (u3);
    \draw (d4) edge[thick, dotted, bend right=20] (u4);
    \draw (d5) edge[thick, dotted, bend right=40] (u5);

    \draw (5, -1.25) node[cop] () {$~$};
    \draw (4.5, -0.75) node[cop] () {$~$};
    \draw (5, -0.25) node[cop] () {$~$};
    \draw (5.5, 0.25) node[cop] () {$~$};

    \draw (5, 2.25) node[prod] () {$~$};
    \draw (3.75, 1) node[prod] () {$~$};
    \draw (6.25, 1) node[prod] () {$~$};
    \draw (5.75, 1.5) node[prod] () {$~$};

    \draw (5, -1.25) -- (3.25, 0.5);
    \draw (5, -1.25) -- (6.75, 0.5);
    \draw (4.5, -0.75) -- (6.25, 1);
    \draw (5, -0.25) -- (3.75, 1);
    \draw (5.5, 0.25) -- (5, 0.75) -- (5.75, 1.5);
    \draw (6.75, 0.5) -- (5, 2.25);
    \draw (3.25, 0.5) -- (5, 2.25);

    \draw (5, -1.5) -- (5, -1.25);
    \draw (5, 2.25) -- (5, 2.5);

    \draw[dotted] (3, 0.5) -- (7, 0.5);

  \end{tikzpicture}
  \caption{Build the sphere and a PC prograph from two triangulations of the $n$-gon.}~\label{glu_trees_example}
\end{figure}

Now, consider two trees with the same number of nodes. Build their 
associated triangulations of the $n$-gon and set opposite orientations on
each of them (as illustrated in Figure~\ref{glu_trees_example}). 
By deformation, we can imagine these two $n$-gons are two hemispheres 
that can be glued along their equator, the root sides being glued together.
This way, considering the underlying trees, the $n-1$ outputs of the tree
of coproducts end up connected to the $n-1$ inputs of the tree of products.
Call the product of the outgoing root side the North and the coproduct of
the ongoing root side the South, and set them at the top and bottom
by convention, respectively.
This operation produces a PC prograph having its products connected up 
to the North point and its coproducts connected up from the South. 
North and South are also connected by the only upside-down
edge of the equator, on the dark side of the sphere. We obtain a figure similar
to Figure~\ref{duality_exemple}.

An immediate property is that PC prographs having coproducts all connected from
the South and products all connected up to the North are exactly PC prographs
avoiding a pattern\,: the output of a product grafted as input of a coproduct.

\begin{thm}
  The subset of PC prographs avoiding edges of types VII (connecting a
  product output to a coproduct input) endowed with the two classical flip types 
  forms a lattice which is exactly the product lattice of the classical Tamari 
  lattice with itself.
\end{thm}

\subsection{PC prographs whose coproducts show up first}

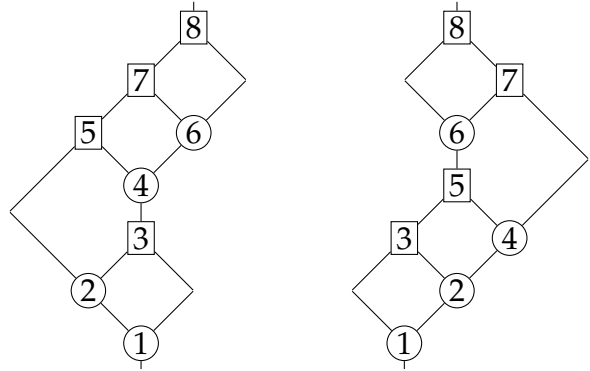
\begin{figure}[h]
  \centering
\begin{tikzpicture}[scale=0.7]
  \tikzset{
    boxi/.style={fill,draw,rectangle,minimum size=5pt,inner sep=1pt}
  }
  \tikzstyle{prod}=[draw,rectangle,minimum size=7pt,inner sep=2pt]
  \tikzstyle{cop}=[draw,circle,minimum size=6pt,inner sep=1pt]
  \tikzstyle{virt}=[minimum size=0pt,inner sep=0pt]  

  \draw (0,0) node[cop] (c1) {$1$};
  \draw (-1,1) node[cop] (c2) {$2$};
  \draw (0,2) node[prod] (p1) {$3$};
  \draw (0,3) node[cop] (c3) {$4$};  
  \draw (-1,4) node[prod] (p2) {$5$};
  \draw (1,4) node[cop] (c4) {$6$};
  \draw (0,5) node[prod] (p3) {$7$};
  \draw (1,6) node[prod] (p4) {$8$};

  \draw (1,1) node[virt] (e1) {};
  \draw (-2.5,2.5) node[virt] (e2) {};
  \draw (2,5) node[virt] (e3) {};
  \draw (0,-0.5) node[virt] (e4) {};
  \draw (1,6.5) node[virt] (e5) {};
  
  \draw (c1) edge (c2);
  \draw (c1) edge (e1);
  \draw (e1) edge (p1);
  \draw (c2) edge (e2);
  \draw (e2) edge (p2);
  \draw (c2) edge (p1);
  \draw (p1) edge (c3);
  \draw (c3) edge (p2);
  \draw (c3) edge (c4);
  \draw (p3) edge (p4);
  \draw (p2) edge (p3);
  \draw (c4) edge (p3);
  \draw (c4) edge (e3);
  \draw (e3) edge (p4);
  \draw (p4) edge (e5);
  \draw (c1) edge (e4);

  \draw (5,0) node[cop] (c1) {$1$};
  \draw (6,1) node[cop] (c2) {$2$};
  \draw (7,2) node[cop] (c3) {$4$};
  \draw (5,2) node[prod] (p1) {$3$};
  \draw (6,3) node[prod] (p2) {$5$};
  \draw (6,4) node[cop] (c4) {$6$};  
  \draw (7,5) node[prod] (p3) {$7$};
  \draw (6,6) node[prod] (p4) {$8$};

  \draw (4,1) node[virt] (e1) {};
  \draw (8.5,3.5) node[virt] (e2) {};
  \draw (5,5) node[virt] (e3) {};
  \draw (5,-0.5) node[virt] (e4) {};
  \draw (6,6.5) node[virt] (e5) {};

  \draw (c1) edge (c2);
  \draw (c1) edge (e1);
  \draw (e1) edge (p1);
  \draw (c3) edge (e2);
  \draw (e2) edge (p3);
  \draw (c2) edge (p1);
  \draw (p1) edge (p2);
  \draw (c3) edge (p2);
  \draw (p2) edge (c4);
  \draw (p3) edge (p4);
  \draw (c2) edge (c3);
  \draw (c4) edge (p3);
  \draw (c4) edge (e3);
  \draw (e3) edge (p4);
  \draw (p4) edge (e5);
  \draw (c1) edge (e4);
\end{tikzpicture}
\caption{Labelling operators of a prograph and its image by the Schützenberger involution.}~\label{labelboxes}
\end{figure}

In Section 5 of~\cite{borie_fpsac_slc_2017}, we did expose a process 
labelling the operators\,: from bottom to top and from left to right 
for operators whose entries lean on already labelled operators. This 
natural process (called Boriefication in~\cite{cordero:tel-03011831}) is 
recalled in Figure~\ref{labelboxes}.

Among all PC prographs, a small portion produce labels $1, 2, \dots , n$
to coproducts and labels $n+1, n+2, \dots , 2n$ to products. In some sense, the 
depth-left-first transversal of operators enumerate all coproducts before the 
first product can go out.
Regarding standard tableaux, it is equivalent to check that the last entry of the 
first row (indexing the last coproduct) is smaller than the first entry of the last row 
(indexing the first product to have its two entries labelled). This subfamily of 
PC prographs is counted by Sequence A274969 of OEIS~\cite{Sloane}. 
It is also 
a subfamily of that described in Section~\ref{glutree}. Such prographs 
can be obtained as a gluing of two binary trees with a reattachment condition.

Labelling nodes of binary trees $T$ with the depth-left-first search algorithm,
we can refine the enumeration of binary trees by counting leaves to the left of
the last node. Denoting by $size(T)$ the number of nodes in $T$ and by $dg(T)$ 
(disabled grafting sites) the number of leaves to the left of the last node, we 
define $Cat_{i}(q)$ to be the following formal sum over all binary 
trees:
\begin{equation}
  Cat_{i}(q) = \displaystyle\sum_{size(T) = i} q^{dg(T)}.
\end{equation}

$Cat_{i}(q)$ deploys one very famous Catalan triangle. Here are the first 
values.
\begin{small}
\begin{displaymath}
  \begin{array}{l}
  Cat_{0}(q) = 1 \\
  Cat_{1}(q) = 1 \\
  Cat_{2}(q) = 1 + q \\
  Cat_{3}(q) = 1 + 2q + 2q^{2} \\
  Cat_{4}(q) = 1 + 3q + 5q^{2} + 5q^{3} \\
  Cat_{5}(q) = 1 + 4q + 9q^{2} + 14q^{3} + 14q^{4} \\
  Cat_{6}(q) = 1 + 5q + 14q^{2} + 28q^{3} + 42q^{4} + 42q^{5}  \\
  \end{array}
\end{displaymath}
\end{small}

It is known that the coefficients of this triangle satisfy
\begin{displaymath}
  C(n, k) = \binom{n+k}{k} - \binom{n+k}{k-1}.
\end{displaymath}

\begin{prop}
  PC prographs of size $n$ whose operators are labelled by indices
  $1, 2, \dots, n$ for coproducts and indices $n+1, n+2, \dots, 2n$
  for products are counted by pairs of binary trees of size $n$
  with $d_1$ and $d_2$ disabled grafting sites and $d_1+d_2 \leqslant
  n+1$.
  Equivalently, we have
\begin{displaymath}
  \sum_{d_1 + d_2 \leqslant n+1} C(n, d_1)C(n, d_2)  = \binom{3n}{n} - 2\binom{3n}{n-1} + \binom{3n}{n-2}
\end{displaymath}
\end{prop}

This formula is more complicated than the compact combination
of three binomial due to Janis Stipins A274969 of OEIS~\cite{Sloane}.
Anyway, it provides an alternative way to compute OEIS A274969 : 
1, 1, 4, 21, 121, 728, 4488, ... 
Furthermore, it exposes an alternative description of pushall stack words 
of length $3n$~\cite{pierrot_2_stack} as pairs of binary trees. The first 
values are the sum of the first coefficients of the square of $Cat_{i}(q)$.

\begin{small}
\begin{displaymath}
  \begin{array}{l}
  Cat_{0}(q)^2~ mod~ q = 1 =_{|q=1} 1 \\
  Cat_{1}(q)^2~ mod~ q^2 = 1 =_{|q=1} 1 \\
  Cat_{2}(q)^2~ mod~ q^3 = 1 + 2q + q^2 =_{|q=1} 4 \\
  Cat_{3}(q)^2~ mod~ q^4 = 1 + 4q + 8q^{2} +8q^{3} =_{|q=1} 21 \\
  Cat_{4}(q)^2~ mod~ q^5 = 1 + 6q + 19q^{2} + 40q^{3} + 55q^{4} =_{|q=1} 121 \\
  Cat_{5}(q)^2~ mod~ q^6 = 1 + 8q + 34q^{2} + 100q^{3} + 221q^{4} + 364q^{5} =_{|q=1} 728 \\
  Cat_{6}(q)^2~ mod~ q^7 = 1 + 10q + 53q^{2} + 196q^{3} + 560q^{4} + 1288q^{5} + 2380q^{6} =_{|q=1} 4488 \\
  \end{array}
\end{displaymath}
\end{small}

\subsection{A poset that is not a lattice on PC prograph}

\begin{prop}
The set of all PC prographs with the four rotation rules does not form a lattice.
\end{prop}

The generalization of the classical Catalan world stops here. The smallest
pair of uncomparable elements are prographs of size $3$. There are $3$ 
elements higher in the poset than both prographs but no join can be defined 
(there are two uncomparable minima among these $3$ higher elements).

\begin{displaymath}
\begin{array}{|c|c|c|} \hline 
  7 & 8 & 9 \\ \hline 
  2 & 4 & 6 \\ \hline 
  1 & 3 & 5 \\ \hline  
\end{array} \qquad
\begin{array}{|c|c|c|} \hline 
  5 & 8 & 9 \\ \hline 
  2 & 6 & 7 \\ \hline 
  1 & 3 & 4 \\ \hline  
\end{array}
\end{displaymath}

\begin{figure}[h]
    \centering
    \begin{tikzpicture}[xscale=1.2, yscale=2.04, every node/.style={scale=0.75}]
  
      \draw (0, 6) node (T6-0) {$\begin{array}{|c|c|c|} \hline 4 & 7 & 9 \\ \hline 2 & 5 & 8 \\ \hline 1 & 3 & 6 \\ \hline  \end{array}$};

      \draw (-4.4, 5) node (T5--4) {$\begin{array}{|c|c|c|} \hline 6 & 7 & 9 \\ \hline 2 & 4 & 8 \\ \hline 1 & 3 & 5 \\ \hline  \end{array}$};
      \draw (T6-0) edge[->, double] (T5--4);
  
      \draw (-2.8, 5) node (T5--3) {$\begin{array}{|c|c|c|} \hline 4 & 8 & 9 \\ \hline 2 & 5 & 7 \\ \hline 1 & 3 & 6 \\ \hline  \end{array}$};
      \draw (T6-0) edge[->, double] (T5--3);
  
      \draw (-1.2, 5) node (T5--1) {$\begin{array}{|c|c|c|} \hline 4 & 6 & 9 \\ \hline 2 & 5 & 8 \\ \hline 1 & 3 & 7 \\ \hline  \end{array}$};
      \draw (T6-0) edge[->, decorate, decoration={snake}] (T5--1);
  
      \draw (1.2, 5) node (T5-1) {$\begin{array}{|c|c|c|} \hline 3 & 7 & 9 \\ \hline 2 & 5 & 8 \\ \hline 1 & 4 & 6 \\ \hline  \end{array}$};
      \draw (T6-0) edge[->, decorate, decoration={snake}] (T5-1);
  
      \draw (2.8, 5) node (T5-3) {$\begin{array}{|c|c|c|} \hline 4 & 7 & 9 \\ \hline 3 & 5 & 8 \\ \hline 1 & 2 & 6 \\ \hline  \end{array}$};
      \draw (T6-0) edge[->, dotted] (T5-3);
  
      \draw (4.4, 5) node (T5-4) {$\begin{array}{|c|c|c|} \hline 5 & 7 & 9 \\ \hline 2 & 6 & 8 \\ \hline 1 & 3 & 4 \\ \hline  \end{array}$};
      \draw (T6-0) edge[->, dotted] (T5-4);

      \draw (-5.5, 4) node (T4--5) {$\begin{array}{|c|c|c|} \hline 6 & 8 & 9 \\ \hline 2 & 4 & 7 \\ \hline 1 & 3 & 5 \\ \hline  \end{array}$};
      \draw (T5--4) edge[->, double] (T4--5);

      \draw (-3.6, 4) node (T4--4) {$\begin{array}{|c|c|c|} \hline 6 & 7 & 9 \\ \hline 3 & 4 & 8 \\ \hline 1 & 2 & 5 \\ \hline  \end{array}$};
      \draw (T5--4) edge[->, dotted] (T4--4);
      \draw (T5-3) edge[->, double] (T4--4);
  
      \draw (-2.4, 4) node (T4--3) {$\begin{array}{|c|c|c|} \hline 5 & 7 & 9 \\ \hline 2 & 4 & 8 \\ \hline 1 & 3 & 6 \\ \hline  \end{array}$};
      \draw (T5--4) edge[->, decorate, decoration={snake}] (T4--3);
  
      \draw (0, 4) node (T4-0) {$\begin{array}{|c|c|c|} \hline 6 & 7 & 9 \\ \hline 2 & 5 & 8 \\ \hline 1 & 3 & 4 \\ \hline  \end{array}$};
      \draw (T5--4) edge[->, dotted] (T4-0);
      \draw (T5-4) edge[->, double] (T4-0);

      \draw (2.4, 4) node (T4-3) {$\begin{array}{|c|c|c|} \hline 4 & 7 & 9 \\ \hline 2 & 6 & 8 \\ \hline 1 & 3 & 5 \\ \hline  \end{array}$};
      \draw (T5-4) edge[->, decorate, decoration={snake}] (T4-3);
      \draw (T4-3) edge[->, decorate, decoration={snake}] (T4--4);
  
      \draw (3.6, 4) node (T4-4) {$\begin{array}{|c|c|c|} \hline 5 & 8 & 9 \\ \hline 2 & 6 & 7 \\ \hline 1 & 3 & 4 \\ \hline  \end{array}$};
      \draw (T5-4) edge[->, double] (T4-4);
      \draw (T5--3) edge[->, dotted] (T4-4);
      \draw (T4--3) edge[->, decorate, decoration={snake}] (T4-4);

      \draw (5.5, 4) node (T4-5) {$\begin{array}{|c|c|c|} \hline 5 & 7 & 9 \\ \hline 3 & 6 & 8 \\ \hline 1 & 2 & 4 \\ \hline  \end{array}$};
      \draw (T5-4) edge[->, dotted] (T4-5);
  
      \draw (-6, 3) node (T3--6) {$\begin{array}{|c|c|c|} \hline 7 & 8 & 9 \\ \hline 2 & 4 & 6 \\ \hline 1 & 3 & 5 \\ \hline  \end{array}$};
      \draw (T4--5) edge[->, double] (T3--6);
      \draw (T5--3) edge[->, double] (T3--6);
  
      \draw (-4.8, 3) node (T3--5) {$\begin{array}{|c|c|c|} \hline 6 & 8 & 9 \\ \hline 2 & 5 & 7 \\ \hline 1 & 3 & 4 \\ \hline  \end{array}$};
      \draw (T4--5) edge[->, dotted] (T3--5);
      \draw (T4-0) edge[->, double] (T3--5);
  
      \draw (-3.6, 3) node (T3--4) {$\begin{array}{|c|c|c|} \hline 5 & 8 & 9 \\ \hline 2 & 4 & 7 \\ \hline 1 & 3 & 6 \\ \hline  \end{array}$};
      \draw (T4--5) edge[->, decorate, decoration={snake}] (T3--4);
      \draw (T4--3) edge[->, double] (T3--4);
  
      \draw (-2.4, 3) node (T3--3) {$\begin{array}{|c|c|c|} \hline 5 & 7 & 9 \\ \hline 3 & 4 & 8 \\ \hline 1 & 2 & 6 \\ \hline  \end{array}$};
      \draw (T4--4) edge[->, decorate, decoration={snake}] (T3--3);
      \draw (T4--3) edge[->, dotted] (T3--3);
      \draw (T5-1) edge[->, decorate, decoration={snake}] (T3--3);
  
      \draw (-1.2, 3.32) node (T3--1) {$\begin{array}{|c|c|c|} \hline 3 & 8 & 9 \\ \hline 2 & 5 & 7 \\ \hline 1 & 4 & 6 \\ \hline  \end{array}$};
      \draw (T5--3) edge[->, decorate, decoration={snake}] (T3--1);
      \draw (T5-1) edge[->, double] (T3--1);
      \draw (-1.2, 2.68) node (T3--1d) {$\begin{array}{|c|c|c|} \hline 5 & 6 & 9 \\ \hline 2 & 4 & 8 \\ \hline 1 & 3 & 7 \\ \hline  \end{array}$};
      \draw (T5--1) edge[->, double] (T3--1d);
      \draw (T4--3) edge[->, decorate, decoration={snake}] (T3--1d);
  
      \draw (0, 3.32) node (T3-0) {$\begin{array}{|c|c|c|} \hline 3 & 6 & 9 \\ \hline 2 & 5 & 8 \\ \hline 1 & 4 & 7 \\ \hline  \end{array}$};
      \draw (T5--1) edge[->, decorate, decoration={snake}] (T3-0);
      \draw (T5-1) edge[->, decorate, decoration={snake}] (T3-0);
      \draw (0, 2.68) node (T3-0d) {$\begin{array}{|c|c|c|} \hline 4 & 8 & 9 \\ \hline 3 & 5 & 7 \\ \hline 1 & 2 & 6 \\ \hline  \end{array}$};
      \draw (T5--3) edge[->, dotted] (T3-0d);
      \draw (T5-3) edge[->, double] (T3-0d);
  
      \draw (1.2, 3.32) node (T3-1) {$\begin{array}{|c|c|c|} \hline 4 & 6 & 9 \\ \hline 3 & 5 & 8 \\ \hline 1 & 2 & 7 \\ \hline  \end{array}$};
      \draw (T5-3) edge[->, decorate, decoration={snake}] (T3-1);
      \draw (T5--1) edge[->, dotted] (T3-1);
      \draw (1.2, 2.68) node (T3-1d) {$\begin{array}{|c|c|c|} \hline 3 & 7 & 9 \\ \hline 2 & 6 & 8 \\ \hline 1 & 4 & 5 \\ \hline  \end{array}$};
      \draw (T5-1) edge[->, dotted] (T3-1d);
      \draw (T4-3) edge[->, decorate, decoration={snake}] (T3-1d);
  
      \draw (2.4, 3) node (T3-3) {$\begin{array}{|c|c|c|} \hline 4 & 8 & 9 \\ \hline 2 & 6 & 7 \\ \hline 1 & 3 & 5 \\ \hline  \end{array}$};
      \draw (T4-4) edge[->, decorate, decoration={snake}] (T3-3);
      \draw (T4-3) edge[->, double] (T3-3);
      \draw (T5--1) edge[->, decorate, decoration={snake}] (T3-3);
  
      \draw (3.6, 3) node (T3-4) {$\begin{array}{|c|c|c|} \hline 4 & 7 & 9 \\ \hline 3 & 6 & 8 \\ \hline 1 & 2 & 5 \\ \hline  \end{array}$};
      \draw (T4-5) edge[->, decorate, decoration={snake}] (T3-4);
      \draw (T4-3) edge[->, dotted] (T3-4);
  
      \draw (4.8, 3) node (T3-5) {$\begin{array}{|c|c|c|} \hline 6 & 7 & 9 \\ \hline 3 & 5 & 8 \\ \hline 1 & 2 & 4 \\ \hline  \end{array}$};
      \draw (T4-5) edge[->, double] (T3-5);
      \draw (T4-0) edge[->, dotted] (T3-5);
  
      \draw (6, 3) node (T3-6) {$\begin{array}{|c|c|c|} \hline 5 & 7 & 9 \\ \hline 4 & 6 & 8 \\ \hline 1 & 2 & 3 \\ \hline  \end{array}$};
      \draw (T5-3) edge[->, dotted] (T3-6);
      \draw (T4-5) edge[->, dotted] (T3-6);
  
      \draw (-5.5, 2) node (T2--5) {$\begin{array}{|c|c|c|} \hline 7 & 8 & 9 \\ \hline 2 & 5 & 6 \\ \hline 1 & 3 & 4 \\ \hline \end{array}$};
      \draw (T3--6) edge[->, dotted] (T2--5);
      \draw (T3--5) edge[->, double] (T2--5);
      \draw (T3--4) edge[->, decorate, decoration={snake}] (T2--5);
      \draw (T4-4) edge[->, double] (T2--5);
  
      \draw (-3.6, 2) node (T2--4) {$\begin{array}{|c|c|c|} \hline 5 & 8 & 9 \\ \hline 3 & 6 & 7 \\ \hline 1 & 2 & 4 \\ \hline  \end{array}$};
      \draw (T4-5) edge[->, double] (T2--4);
      \draw (T4-4) edge[->, dotted] (T2--4);
      \draw (T3--3) edge[->, decorate, decoration={snake}] (T2--4);
  
      \draw (-2.4, 2) node (T2--3) {$\begin{array}{|c|c|c|} \hline 5 & 8 & 9 \\ \hline 3 & 4 & 7 \\ \hline 1 & 2 & 6 \\ \hline  \end{array}$};
      \draw (T3--3) edge[->, double] (T2--3);
      \draw (T3--4) edge[->, dotted] (T2--3);
  
      \draw (0, 2) node (T2-0) {$\begin{array}{|c|c|c|} \hline 6 & 8 & 9 \\ \hline 3 & 5 & 7 \\ \hline 1 & 2 & 4 \\ \hline  \end{array}$};
      \draw (T3-5) edge[->, double] (T2-0);
      \draw (T3--5) edge[->, dotted] (T2-0);
  
      \draw (2.4, 2) node (T2-3) {$\begin{array}{|c|c|c|} \hline 4 & 8 & 9 \\ \hline 3 & 6 & 7 \\ \hline 1 & 2 & 5 \\ \hline  \end{array}$};
      \draw (T2--4) edge[->, decorate, decoration={snake}] (T2-3);
      \draw (T3-3) edge[->, dotted] (T2-3);
      \draw (T3-4) edge[->, double] (T2-3);
  
      \draw (3.6, 2) node (T2-4) {$\begin{array}{|c|c|c|} \hline 6 & 8 & 9 \\ \hline 3 & 4 & 7 \\ \hline 1 & 2 & 5 \\ \hline  \end{array}$};
      \draw (T4--5) edge[->, dotted] (T2-4);
      \draw (T4--4) edge[->, double] (T2-4);
      \draw (T3-3) edge[->, decorate, decoration={snake}] (T2-4);
      \draw (T2-4) edge[->, decorate, decoration={snake}] (T2--3);
  
      \draw (5.5, 2) node (T2-5) {$\begin{array}{|c|c|c|} \hline 6 & 7 & 9 \\ \hline 4 & 5 & 8 \\ \hline 1 & 2 & 3 \\ \hline \end{array}$};
      \draw (T3-6) edge[->, double] (T2-5);
      \draw (T3-5) edge[->, dotted] (T2-5);
      \draw (T3-4) edge[->, decorate, decoration={snake}] (T2-5);
      \draw (T4--4) edge[->, dotted] (T2-5);
  
      \draw (-4.4, 1) node (T1--4) {$\begin{array}{|c|c|c|} \hline 7 & 8 & 9 \\ \hline 3 & 5 & 6 \\ \hline 1 & 2 & 4 \\ \hline  \end{array}$};
      \draw (T2--5) edge[->, dotted] (T1--4);
      \draw (T2-0) edge[->, double] (T1--4);
      \draw (T2--4) edge[->, double] (T1--4);
      \draw (T2--3) edge[->, decorate, decoration={snake}] (T1--4);
  
      \draw (-2.8, 1) node (T1--3) {$\begin{array}{|c|c|c|} \hline 7 & 8 & 9 \\ \hline 3 & 4 & 6 \\ \hline 1 & 2 & 5 \\ \hline  \end{array}$};
      \draw (T3--6) edge[->, dotted] (T1--3);
      \draw (T3-0d) edge[->, double] (T1--3);
      \draw (T3--1d) edge[->, decorate, decoration={snake}] (T1--3);
      \draw (T2-4) edge[->, double] (T1--3);
  
      \draw (-1.2, 1) node (T1--1) {$\begin{array}{|c|c|c|} \hline 3 & 8 & 9 \\ \hline 2 & 6 & 7 \\ \hline 1 & 4 & 5 \\ \hline  \end{array}$};
      \draw (T3-0) edge[->, decorate, decoration={snake}] (T1--1);
      \draw (T3--1) edge[->, dotted] (T1--1);
      \draw (T3-1d) edge[->, double] (T1--1);
      \draw (T3-3) edge[->, decorate, decoration={snake}] (T1--1);
  
      \draw (1.2, 1) node (T1-1) {$\begin{array}{|c|c|c|} \hline 5 & 6 & 9 \\ \hline 3 & 4 & 8 \\ \hline 1 & 2 & 7 \\ \hline  \end{array}$};
      \draw (T3-0) edge[->, decorate, decoration={snake}] (T1-1);
      \draw (T3-1) edge[->, double] (T1-1);
      \draw (T3--1d) edge[->, dotted] (T1-1);
      \draw (T3--3) edge[->, decorate, decoration={snake}] (T1-1);
  
      \draw (2.8, 1) node (T1-3) {$\begin{array}{|c|c|c|} \hline 5 & 8 & 9 \\ \hline 4 & 6 & 7 \\ \hline 1 & 2 & 3 \\ \hline  \end{array}$};
      \draw (T3-6) edge[->, double] (T1-3);
      \draw (T3-0d) edge[->, dotted] (T1-3);
      \draw (T3-1d) edge[->, decorate, decoration={snake}] (T1-3);
      \draw (T2--4) edge[->, dotted] (T1-3);
  
      \draw (4.4, 1) node (T1-4) {$\begin{array}{|c|c|c|} \hline 6 & 8 & 9 \\ \hline 4 & 5 & 7 \\ \hline 1 & 2 & 3 \\ \hline  \end{array}$};
      \draw (T2-5) edge[->, double] (T1-4);
      \draw (T2-0) edge[->, dotted] (T1-4);
      \draw (T2-4) edge[->, dotted] (T1-4);
      \draw (T2-3) edge[->, decorate, decoration={snake}] (T1-4);
  
      \draw (0, 0) node (T0-0) {$\begin{array}{|c|c|c|} \hline 7 & 8 & 9 \\ \hline 4 & 5 & 6 \\ \hline 1 & 2 & 3 \\ \hline  \end{array}$};
      \draw (T1--1) edge[->, decorate, decoration={snake}] (T0-0);
      \draw (T1-1) edge[->, decorate, decoration={snake}] (T0-0);
      \draw (T1--3) edge[->, dotted] (T0-0);
      \draw (T1-3) edge[->, double] (T0-0);
      \draw (T1--4) edge[->, dotted] (T0-0);
      \draw (T1-4) edge[->, double] (T0-0);
  
    \end{tikzpicture}
    \caption{The poset formed by the $42$ PC prographs of size $3$.}~\label{assos3tab}
  \end{figure}
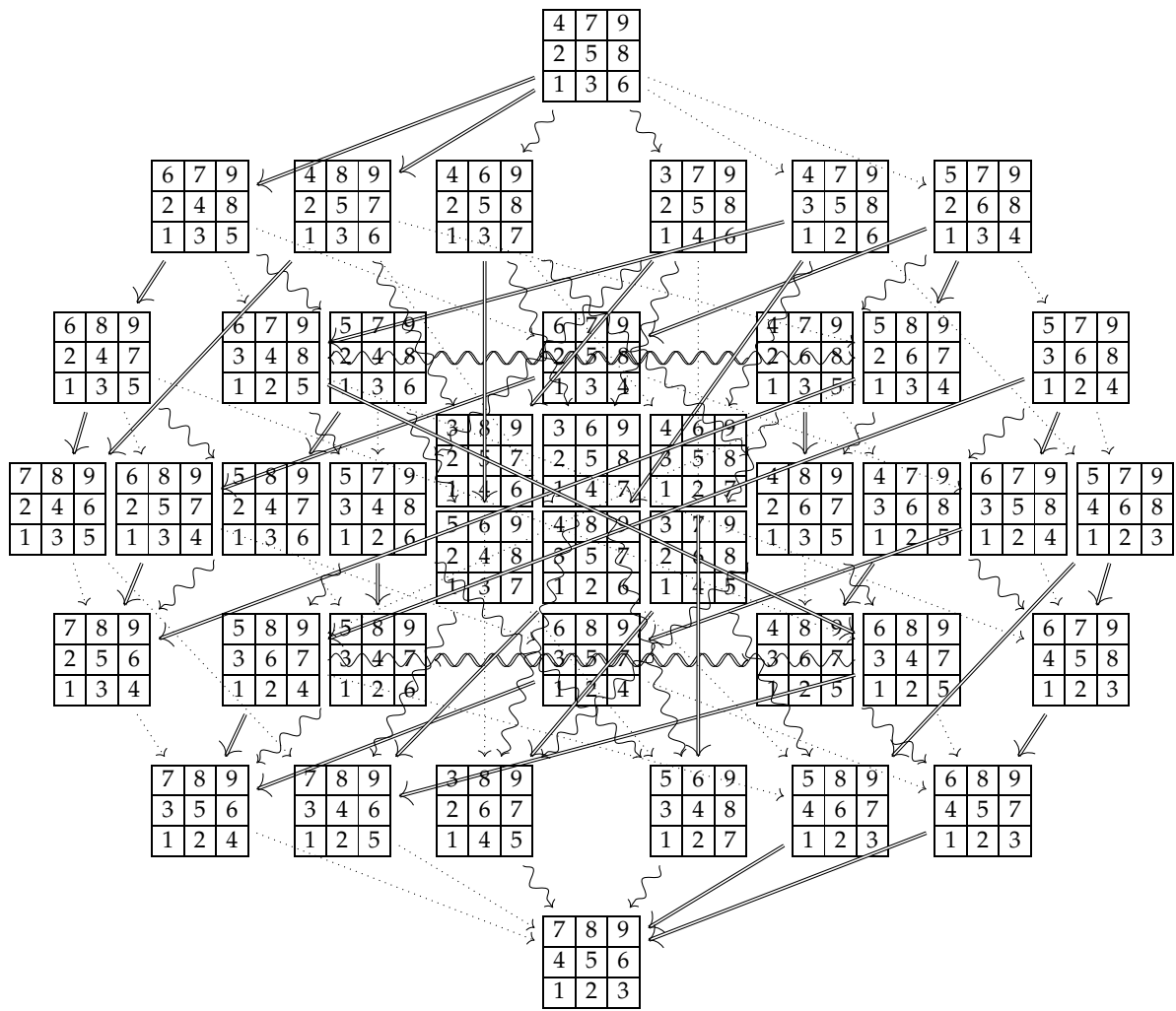

Figure~\ref{assos3tab} displays the poset over the $42$ PC prographs 
of size $3$ (as standard tableaux).



This research was driven by computer exploration using the open-source
mathematical software \texttt{Sage}~\cite{sage} and its algebraic
combinatorics features developed by the \texttt{Sage-Combinat}
community~\cite{Sage-Combinat}.

\printbibliography

@Misc{Sage-Combinat,
      Author = {The {S}age-{C}ombinat community},
      Title = {{S}age-{C}ombinat: enhancing {S}age as a toolbox for computer exploration in algebraic combinatorics},
      note= {{\tt http://combinat.sagemath.org}},
      Year = 2008}

@manual{sage,
  Key          = {Sage},
  Author       = {The Sage Developers},
  Title        = {{S}age {M}athematics {S}oftware ({V}ersion 6.4)},
  note         = {{\tt http://www.sagemath.org}},
  Year         = {2014},
}

@article {Sloane,
     AUTHOR = {Sloane, N. J. A.},
      TITLE = {The on-line encyclopedia of integer sequences},
    JOURNAL = {Notices Amer. Math. Soc.},
   FJOURNAL = {Notices of the American Mathematical Society},
     VOLUME = {50},
       YEAR = {2003},
     NUMBER = {8},
      PAGES = {912--915},
       ISSN = {0002-9920},
      CODEN = {AMNOAN},
    MRCLASS = {11Y55 (05-00 11-00 11B83)},
   MRNUMBER = {1992789 (2004f:11151)},
 MRREVIEWER = {Christian Krattenthaler},
}

@Article{borie_fpsac_slc_2017,
 Author = {Nicolas {Borie}},
 Title = {{Three-dimensional Catalan numbers and product-coproduct prographs.}},
 FJournal = {{S\'eminaire Lotharingien de Combinatoire}},
 Journal = {{S\'emin. Lothar. Comb.}},
 ISSN = {1286-4889/e},
 Volume = {78B},
 Pages = {78b.39, 12},
 Year = {2017},
 Publisher = {Universit\"at Wien, Fakult\"at f\"ur Mathematik, Wien},
 Language = {English},
 MSC2010 = {05E10},
 Zbl = {1384.05165}
}

@unpublished{pierrot_2_stack,
  TITLE = {{2-stack pushall sortable permutations}},
  AUTHOR = {Pierrot, Adeline and Rossin, Dominique},
  URL = {https://hal.archives-ouvertes.fr/hal-00801861},
  NOTE = {41 pages},
  YEAR = {2013},
  MONTH = Mar,
  HAL_ID = {hal-00801861},
  HAL_VERSION = {v1},
}

@article {bipolar,
    AUTHOR = {Kenyon, Richard and Miller, Jason and Sheffield, Scott and
              Wilson, David B.},
     TITLE = {Bipolar orientations on planar maps and {${\rm SLE}_{12}$}},
   JOURNAL = {Ann. Probab.},
  FJOURNAL = {The Annals of Probability},
    VOLUME = {47},
      YEAR = {2019},
    NUMBER = {3},
     PAGES = {1240--1269},
      ISSN = {0091-1798},
   MRCLASS = {60J67 (05C10 28C20 82B20)},
  MRNUMBER = {3945746},
MRREVIEWER = {Nam-Gyu Kang},
       DOI = {10.1214/18-AOP1282},
       URL = {https://doi.org/10.1214/18-AOP1282},
}

@InCollection{MBM_planar_maps,
 Author = {Mireille {Bousquet-M\'elou}},
 Title = {{Counting planar maps, coloured or uncoloured.}},
 BookTitle = {{Surveys in combinatorics 2011. Papers from the 23rd British combinatorial conference, Exeter, UK, July 3--8, 2011}},
 ISBN = {978-1-107-60109-3/pbk},
 Pages = {1--49},
 Year = {2011},
 Publisher = {Cambridge: Cambridge University Press},
 Language = {English},
 MSC2010 = {05C30 05C10 05C15 05A15},
 Zbl = {1244.05118}
}

@phdthesis{cordero:tel-03011831,
  TITLE = {{Explorations combinatoires des structures arborescentes et libres}},
  AUTHOR = {Cordero, Christophe},
  URL = {https://tel.archives-ouvertes.fr/tel-03011831},
  NUMBER = {2019PESC2046},
  SCHOOL = {{Universit{\'e} Paris-Est}},
  YEAR = {2019},
  MONTH = Dec,
  TYPE = {Theses},
  HAL_ID = {tel-03011831},
  HAL_VERSION = {v1},
}

@article {falque_hard_bij,
    AUTHOR = {Falque, Justine},
     TITLE = {A bijection between weighted {D}yck paths and 1234-avoiding
              alternating permutations},
   JOURNAL = {S\'{e}m. Lothar. Combin.},
  FJOURNAL = {S\'{e}minaire Lotharingien de Combinatoire},
    VOLUME = {85B},
      YEAR = {2021},
     PAGES = {Art. 71, 13},
   MRCLASS = {05A19 (05A05 05A15)},
  MRNUMBER = {4311952},
}

\end{document}